\documentclass{aptpub}

  \usepackage{xcolor}
\usepackage{graphicx}
\usepackage[margin=1in]{geometry}
\usepackage[linesnumbered,ruled]{algorithm2e}

\usepackage{empheq}
\usepackage[sort,compress,numbers]{natbib}

\usepackage{subfigure}

  \newtheorem{coro}{Corollary}

\newcommand{\tkappa}{{\tilde{\kappa}}}
\newcommand{\tep}{{\tilde{\epsilon}}}

\authornames{Li and Xu} 
\shorttitle{Uniformly Efficient Simulation for Gaussian  Fields} 

\begin{document}

\title{Uniformly Efficient Simulation for Extremes of Gaussian Random Fields} 

\authorone[University of Minnesota]{Xiaoou Li} 
\authortwo[University of Michigan]{Gongjun Xu}
\addressone{School of Statistics, 224 Church ST SE, Minneapolis, MN, USA, 55455. Email:  lixx1766@umn.edu} 
\addresstwo{Department of Statistics, 1085 South University, Ann Arbor, MI, USA, 48109. Email: gongjun@umich.edu}
\begin{abstract}
This paper considers the problem of simultaneously estimating  rare-event probabilities for  a class of Gaussian random fields.  
 A conventional rare-event simulation method is usually  tailored   to  a specific  rare event  and consequently would lose  estimation efficiency  for  different  events of interest, which  often  results in additional computational cost in such   simultaneous estimation problem. 
To overcome this issue, we propose  a uniformly efficient estimator for a general family of H\"older continuous Gaussian random fields.
We establish the asymptotic and uniform efficiency  of the proposed method  and also conduct simulation studies  to illustrate its effectiveness. 
\end{abstract}



\section{Introduction}
\label{sec:intro}

Consider a continuous Gaussian random field  $\{f(t): t\in T\}$ with  zero mean and unit variance, living on a $d$-dimensional compact set $T\subset R^d$; that is, for
every finite subset of $\{t_1,...,t_n\}\subset T$,
$(f(t_1),...,f(t_n))$ is a multivariate Gaussian random
vector with $Ef(t_i)=0$ and $Var(f(t_i))=1$ for $i=1,\cdots,n$. 
We are interested in estimating the tail probability
$$w_{\sigma,\mu}(b)=P\left(\sup_{t\in T} \{\sigma(t)f(t) +\mu(t)\}>b\right), \mbox{ as } b\rightarrow \infty,$$
simultaneously for a class of continuous mean and variance functions $\mu(t)$ and  $\sigma^2(t)$, where the functions $\mu(t)$ and $\sigma^2(t)$ may be unspecified and only known to be in certain ranges.

The extremes of Gaussian random fields have wide applications in finance, spatial analysis, physical oceanography, and many other disciplines \cite{adler1996stochastic,adler2010applications}.
Tail probabilities of the extremes have been extensively studied in the literature, with its focus mostly on the development of approximations and bounds for the suprema 
\cite[e.g.,][]{LS70,MS70,ST74,CIS,Bor75,Bor03,LT91,TA96,Berman85,Sun93,AW08,AdlTay07,cheng2012mean,debicki2014extremes}. 
 Tail probabilities of  other convex functions of Gaussian random fields have also been studied; see \cite{Liu10,LiuXu11,LiuXu12,LiLiuXuMathOR}.

Most of the sharp theoretical approximations developed in the literature 
  require the evaluation of certain constants that are hard to estimate, such as the Lipschitz-Killing curvatures and Pickands' constant. 
Moreover, although the asymptotic results may provide good approximations for
large tail values as $b\to\infty$, 
evaluation of the approximation results for finite $b$ may be challenging  and it is often unclear how large the tail
values are required to ensure the approximations within an
acceptable range relative to the quantity of interest. 
Therefore, to evaluate the tail probabilities, rare-event simulation serves as an appealing alternative from a computational point of view.
In particular, the design and the analysis do not require very sharp approximations of the tail probabilities.
Importance sampling based efficient simulation procedures have been proposed in the literature to estimate the tail probabilities.
Numerical methods for rare-event analysis of the suprema were studied in \cite{ABL08,ABL09}; see also \cite{azais2009level,li2013rare,LiuXu13,XuWSC14,liuxu2012rare,LiuXu12,LiuXuTomacs} for related studies.

To design asymptotically efficient importance sampling estimator, one needs to construct a change of measure that is tailored to a specific event. Such  construction usually  requires detailed information of the Gaussian random fields, such as $\mu(t)$ and $\sigma(t)$ whose computations themselves are sometimes intensive. In addition, the specific form of the change of measure is sensitive to $\mu(t)$ and $\sigma(t)$ in the sense that the entire simulation needs to be redone even if there is a tiny change of the system. 
This often leads to additional computational overhead especially at the exploratory stage when one often needs to tune different model parameters.
This motivates us to seek for a single Monte Carlo scheme that is efficient for a class of distributions.
An advantage of such uniformly efficient methods is that there is no need to regenerate samples if there is a change in the original system and one just needs to recompute the importance weights. This could save substantial computational time.
Moreover,  this can help researchers efficiently estimate many probabilities for a certain range of mean and variance parameter values, which are often of practical importance.
For instance, in finance risk analysis, there is often uncertainty surrounding the true population values for the mean and variance; portfolio credit risk management may require the estimation of the tail probabilities of extremes for a family of Gaussian processes; in physical system reliability analysis, we may need to evaluate the failure probability for a range of system parameters.

To address the above issues, this study focuses on the problem of simultaneous efficient estimation of $w_{\sigma,\mu}(b)$ for all possible $\mu(t)\in [\mu_l,\mu_u]$ and  $\sigma^2(t)\in [\sigma^2_l, \sigma^2_u]$, $t\in T$, where  $\mu_l\leq\mu_u\in R$ and $\sigma_l\leq\sigma_u\in (0,\infty)$ are constants that are prespecified.
 We propose a mixture type  change of measure that yields uniformly efficient estimation (criterion defined in Section 2). In particular, the uniform efficiency result holds for general H\"older continuous Gaussian random fields and therefore it is applicable to most of the practical problems. 

The remainder of the paper is organized as follows. In Section \ref{notations} we introduce some notions of efficiency and computational complexity under the setting of rare-event simulation. Section \ref{main} provides the construction of our importance sampling estimator and shows the main properties of our algorithm.  Numerical simulations are conducted in Section \ref{Sim} and detailed proofs of our main theorems are given in Section \ref{SecProof}.

\section{Efficiency Criteria}\label{notations}\label{section:efficiency}

\subsection{Efficiency of rare-event simulation and importance sampling}
We first introduce some general notions  of rare-event simulations. Given that the tail probability $w_{\sigma,\mu}(b)$ converges to zero, it is usually meaningful to consider the relative error of a Monte Carlo estimator $L(b)$ with respect to $w_{\sigma,\mu}(b)$. This is because a trivial estimator $L^*(b)\equiv 0$ has an error
$|L^*(b)-w_{\sigma,\mu}(b)|=w_{\sigma,\mu}(b)\rightarrow 0$. In the literature of rare-event simulation  (e.g., \cite{ASMGLY07,juneja2006rare,ABL09}), one usually employs the concept of \emph{polynomially efficiency} as an efficiency criterion.

\begin{definition}[Polynomial efficiency]
An estimator $L(b)$ is said to be \emph{polynomially efficient with the order $q$} in estimating $w_{\sigma,\mu}(b)$ if $EL(b) = w_{\sigma,\mu}(b)$ and there exist  constants $q\geq 0$   and $b_0\geq 0$ such that
\begin{equation}\label{logeff}
  \mbox{sup}_{b\geq b_0} \frac{ Var (L(b))}{|\log w_{\sigma,\mu}(b)|^qw^{2}_{\sigma,\mu}(b)}<\infty.
\end{equation}
When $q=0$, $L(b)$ is also called strongly efficient.
\end{definition}
To illustrate this efficiency criterion, we compare a polynomially efficient estimator with a standard Monte Carlo estimator.
Suppose that we want to estimate $w_{\sigma,\mu}(b)$ with certain relative accuracy with a high probability. That is, we would like to  have an estimator $Z(b)$ such that for some prescribed $\varepsilon, \delta>0$,
\begin{equation}\label{mse}
P\left(\left|Z(b)/ {w_{\sigma,\mu}(b)}-1\right| > \varepsilon\right)< \delta.
\end{equation}
If a  standard Monte Carlo simulation method is used, then it requires at least $n=O(\varepsilon^{-2}\delta ^{-1} \allowbreak w_{\sigma,\mu}^{-1}(b))$ i.i.d. replicates, according to the central limit theorem.
By the Borell-TIS lemma (Lemma \ref{LemBorell}), we know $w_{\sigma,\mu}(b)\leq\exp\{-(1+o(1))b^2/(2\sup_{t\in T}\sigma^2(t))\} $. Therefore, $n$ has to grow at an exponential rate in $b^2$. On the contrary, 
suppose that a polynomially efficient estimator of $w_{\sigma,\mu}(b)$ has been obtained, denoted by $L(b)$. Let $\{L^{(j)}(b): j=1,...,n\}$ be $n$ i.i.d.\ copies of $L(b)$.  Then the averaged estimator
$Z(b) = \frac 1 n \sum_{j=1}^{n}L^{(j)}(b)$
has a mean squared error (MSE)
$
E (Z(b)-w_{\sigma,\mu}(b))^2=
Var (L(b))/n.$ A direct application of Chebyshev's inequality yields
\begin{equation}\label{eq:Cheby}
	P(|Z(b)/w_{\sigma,\mu}(b) -1 |\geq \varepsilon)\leq \frac{Var(L(b))}{n\varepsilon^{2} w_{\sigma,\mu}^{2}(b)}.
\end{equation}
Thus, if $L(b)$ is  a polynomially efficient estimator with the order $q$, it suffices to simulate $n=\varepsilon^{-2}\delta ^{-1}|\log w_{\sigma,\mu}(b)|^q =O(\varepsilon^{-2}\delta ^{-1} b^{2q})$ 
i.i.d. replicates  of $L(b)$ to achieve the accuracy in \eqref{mse}.
Compared with the standard Monte Carlo simulation, polynomially efficient estimators  reduce the computational cost substantially for large $b$.
\begin{remark}\label{remark:weak}
	In the rare event analysis literature, another widely used efficiency criterion is the \emph{weakly efficient} (\cite{ASMGLY07}). An estimator $L(b)$ is said to be weakly efficient in estimating $w_{\sigma,\mu}(b)$, if $EL(b)=w_{\sigma,\mu}(b)$ and for all positive constants $\varepsilon>0$, 
	\begin{equation*}
		\limsup_{b\to\infty} \frac{Var(L(b))}{w^{2-\varepsilon}_{\sigma,\mu}(b)}=0.
	\end{equation*}
	It is easy to verify that if $L(b)$ is polynomially efficient, then $L(b)$ is also weakly efficient. That is, polynomial efficiency is a stronger criterion than the weak efficiency.
\end{remark}

To construct polynomially efficient estimators,  importance sampling is a commonly used method for the variance reduction.
In particular, we have
$$w_{\sigma,\mu}(b)= E \Big[I{ \Big(\sup_{t\in T} \{\sigma(t)f(t) +\mu(t)\}>b \Big)}\Big]
=E^Q \Big[\frac{dP}{dQ}I{ \Big(\sup_{t\in T} \{\sigma(t)f(t) +\mu(t)\}>b \Big)}\Big],$$
where $I{(\cdot)}$ denotes the indicator function, $Q$ is a probability measure that is absolutely continuous with respect to $P$ on the set $\{\sup_{t\in T} \{\sigma(t)f(t) +\mu(t)\}>b\}$, and we use  $E$ and $E^Q$ to denote  the expectations under the measures $P$ and $Q$, respectively.
Then, the random variable defined by
\begin{equation}\label{estr}
L_{\sigma,\mu}(b)=\frac{dP}{dQ}I{ \Big(\sup_{t\in T} \{\sigma(t)f(t) +\mu(t)\}>b \Big)}
\end{equation}
is an unbiased estimator of $w_{\sigma,\mu}(b)$ under the measure $Q$.
To have an efficient estimator, we want to choose $Q$ such that the variance $Var^Q(L_{\sigma,\mu}(b))$ is small. It is straightforward to show that the optimal change of measure is the conditional probability $Q^*(\cdot):=P(\cdot \mid\sup_{t\in T} \{\sigma(t)f(t) +\mu(t)\}>b)=P(\cdot\cap \{\sup_{t\in T} \{\sigma(t)f(t) +\mu(t)\}>b\})/w_{\sigma,\mu}(b)$, for which the corresponding importance sampling estimator has a zero variance. However, $Q^*$ cannot be implemented in practice because  $w_{\sigma,\mu}(b)$, the probability of interest, is unknown beforehand. Therefore, constructing an efficient change of measure usually involves analysis and approximation of the optimal change of measure $Q^*$.
\label{page:zero-variance} 
\subsection{Non-uniformly efficient issue and an example} 
\label{sec:non-u}
Various importance sampling estimators for rare-event analysis of the suprema of Gaussian random fields have been  studied in \cite{ABL08,ABL09,azais2009level,li2013rare}.
As the measure $Q^*$ depends on the mean and variance function $\sigma(\cdot)$ and $\mu(\cdot)$, the designed measures  usually depend on the $\mu(\cdot)$ and $\sigma(\cdot)$ as well.
As a consequence, a measure $Q$ that gives  an efficient estimator $L_{\sigma,\mu}(b)=\frac{dP}{dQ}I{(\sup_{t\in T} \{\sigma(t)f(t) +\mu(t)\}>b)}$
 for $w_{\sigma,\mu}(b)$ may not be  efficient any more for estimating $w_{\sigma',\mu'}(b)$, where $\sigma'(t)$ and $\mu'(t)$ are two different variance and mean functions. That is, the corresponding importance sampling estimator based on $Q$
 $$L_{\sigma',\mu'}(b):=\frac{dP}{dQ}I{(\sup_{t\in T} \{\sigma'(t)f(t) +\mu'(t)\}>b)}$$
 may not be an efficient estimator for $w_{\sigma',\mu'}(b)$ .

To illustrate the non-uniform efficiency issue, we take the estimator proposed in  \cite{ABL09}  as an example.
For simplicity, we consider the case when $T$ contains finite points and write $T:= \{t_{1},\cdots, t_M\}$.

	For known $\mu$ and $\sigma$, \cite{ABL09} proposed the following simulation procedure in Algorithm 1.
\begin{algorithm}[ht]
\caption{Sampling procedure proposed by \cite{ABL09}}
\label{alg:nonu}
\KwIn{$T= \{t_{1},\cdots, t_M\}$.}
Simulate a random variable $\tau \in\{t_{1},\cdots, t_M\}$ according to the following probability measure:
 \begin{equation}\label{taub}
 P(\tau=t_i) = \frac{P(\sigma(t_i) f(t_i)+\mu(t_i)>b)}{\sum_{j=1}^M
 P(\sigma(t_j) f(t_j)+\mu(t_j)>b)};
 \end{equation}

Given the realized $\tau$, simulate $f(\tau)$ conditional on $\sigma(\tau) f(\tau)+\mu(\tau)>b$\;\label{line2}
Given $(\tau,f(\tau))$, simulate the rest  $\{f(t): t\neq \tau, t\in T\}$ from the original conditional distribution under $P$.

\KwOut{$f(t)$ for $t\in T$}
 \end{algorithm}
Let $Q^{\dagger}$ be the corresponding change of measure. We have
$$\frac{dQ^{\dagger}}{dP}=\frac{\sum_{i=1}^M I{(\sigma(t_i) f(t_i)+\mu(t_i)>b)}}
{\sum_{i=1}^M P(\sigma(t_i) f(t_i)+\mu(t_i)>b)}.$$
\cite{ABL09} showed that $L_{\sigma,\mu}(b)=\frac{dP}{dQ^{\dagger}}I{(\sup_{t\in T} \{\sigma(t)f(t) +\mu(t)\}>b)}$ is a polynomially efficient estimator for
$w_{\sigma,\mu}(b)$ with the order $q=0$. We explain intuitively why this estimator is efficient. First, Algorithm~\ref{alg:nonu} samples a random index $\tau$ whose distribution is approximating that of $t^*:=\arg\max_{t_i}(\sigma(t_i)f(t_i)+\mu(t_i))$. Second, it simulates $f(\tau)$ approximately from the conditional distribution $P(f(t^*)\in\cdot| f(t^*)>b)$. Third, Algorithm~\ref{alg:nonu} simulates the $f(t)$ at $t\neq \tau$ according to the original conditional distribution given $(f(\tau),\tau)$. Combining the three steps, the entire sample path $\{f(t): t\in T \}$ generated from Algorithm~\ref{alg:nonu} approximately follows the conditional distribution $\{
f(t): t\in T| \max_{t_i}(\sigma(t_i)f(t_i)+\mu(t_i))>b\}$. According to the discussion on page~\pageref{page:zero-variance}, this conditional probability measure is the optimal change of measure. See \cite{ABL09} for rigorous justifications of the above statements.

 Let $\mu'$ and $\sigma'$ be a different mean and variance function. We have  Proposition \ref{prop:nonu} for the estimator $$L_{\sigma',\mu'}(b):=\frac{dP}{dQ^{\dagger}}I{ \Big(\sup_{t\in T} \{\sigma'(t)f(t) +\mu'(t)\}>b \Big)}.$$ 	
\begin{proposition}\label{prop:nonu}
	Let $\mu'(t)=\mu(t)=0$ for all $t\in T$. 
	\begin{itemize}
		\item [(i)] If  $\sigma'(t)\leq \sigma(t)$ for all $t\in T$ and $\max_{t_i\in T}\sigma'(t_i)<\max_{t_i\in T}\sigma(t_i)$, then 
		for some constant $\varepsilon>0$,
$$
\lim_{b\to\infty}\frac{E^{Q^{\dagger}}\Big[\left(\frac{dP}{dQ^{\dagger}}\right)^{2};
{\max}_{t_i\in T}\sigma'(t_i)f(t_i)>b\Big]}
{ w^{2-\epsilon}_{\sigma',\mu}(b)} =
\infty.$$
\item[(ii)] If $\max_{t_i\in T}\sigma'(t_i)>\max_{t_i\in T}\sigma(t_i)$, then $\frac{dP}{dQ^{\dagger}}$ is not well defined on the event $ \{{\max}_{t_i\in T} \sigma'(t_i) f(t_i)>b\}\allowbreak$.
	\end{itemize}
\end{proposition}

 According to the definition of weakly efficient estimator in Remark~\ref{remark:weak}, the first part of the above proposition implies that $L_{\sigma,\mu}(b)$ is not weakly efficient for estimating $w_{\sigma',\mu'}(b)$ if $ \max_{t_i\in T}\sigma'(t_i)>\max_{t_i\in T}\sigma(t_i)$, and is therefore not polynomially efficient.
The second part of the above proposition implies that the estimator $L_{\sigma,\mu}(b)$ is not well defined when $ \max_{t_i\in T}\sigma'(t_i)>\max_{t_i\in T}\sigma(t_i)$. Therefore, for each $L_{\sigma,\mu}(b)$ there always exist mean and variance functions $\mu'(\cdot)$, $\sigma'(\cdot)$ such that $\mu'(t)\in[\mu_l,\mu_u], \sigma'(t)\in[\sigma_l,\sigma_u]$ and $L_{\sigma,\mu}(b)$ is \emph{not} (weakly) efficient for estimating $w_{\sigma',\mu'}(b)$. 
We use a simple numerical study to further illustrate this.
\begin{example}\label{example:iid}
	Consider i.i.d. standard normal random variables $\{f(t), t=1,\cdots,100\}$.
For simplicity, we take $\mu(t)=0$ and $\sigma(t)=\sigma$ for all $t$. The probability of interest is $P(\sigma\max_t f(t) >b)$ for $\sigma\in[0.3,1] $ and $b=3$. This is equivalent to simulating $P(\max_t f(t) >b)$ for all $b\in[3,10]$. 
Table \ref{result2}  displays the simulation results for
$\sigma = 0.3, 0.6$  and $1$, from Algorithm~\ref{alg:nonu}, where the change of measure is constructed based on $\sigma=1$. The results are based on $10^4$ independent simulations.
 We report the estimated tail probability (est.), the estimated standard deviation (sd.) of $L_{\sigma,\mu}(b)$, and  the coefficient of variation (CV), which is the ratio sd./est.. 
We also give the theoretical values of the tail probabilities, that is, $P(\max_i f(t_i) >b/\sigma) = 1- {\Phi}(b/\sigma)^{100}$ where ${\Phi}(x)=\int_{-\infty}^x \frac{1}{\sqrt{2\pi}}e^{-t^2/2}dt$ denotes the left tail probability of the standard Gaussian distribution. 
We can see that the estimator is more efficient when $\sigma$ value is equal to the designed value 1 and less  for other $\sigma$ values. In particular, when $\sigma=0.3$, it gives 0 estimated value.

\begin{table}[h]
\centering
\begin{tabular}{ccccc}\hline
$\sigma$&est.&sd.&CV & Theoretical Value \\
\hline
 0.3 &  0 & 0 &   NA & 7.62e-22\\
 0.6&  1.35e-05 & 1.35e-03 &  1.00e+02 & 2.87e-05\\
 1 &  1.26e-01 & 2.32e-02 &  1.84e-01 & 1.26e-01\\
 \hline
\end{tabular}
\caption{Estimates based on Algorithm~\ref{alg:nonu}}
\label{result2}
\end{table}
\end{example}
The above non-uniform efficiency result can be extended, with similar techniques, to the importance sampling estimators in  \cite{ABL09}  when $\{f(t), t\in T\}$ is a continuous Gaussian random field. It can also be extended to the case when other change of measures are used such as \cite{li2013rare}. 
In general, if the construction of a rare-event change of measure relies heavily on the mean and variance functions, then it
would not be efficient for another set of functions.

\subsection{Uniform Efficiency} 

In applications, one is often interested in estimating many probabilities for a certain range of mean and variance parameter values, such as evaluating   the tail probabilities of a loss distribution for a range of loss thresholds in portfolio credit risk management
(e.g., \citep{glasserman2013monte,glasserman2008uniformly}).  This motivates us to construct a change of measure such that the corresponding importance sampling estimator $L_{\sigma,\mu}(b)$ is polynomially efficient for a family of functions $\mu$ and $\sigma$. In particular,  this paper  considers $\mu$ and $\sigma$ satisfying the following  condition:

\begin{itemize}
 	\item [C1.]For all $t\in T$, $\mu(t)\in [\mu_l,\mu_u]$  and $\sigma^2(t)\in [\sigma^2_l, \sigma^2_u].$ Moreover, $\mu$ and $\sigma$ are H\"older continuous in the sense that there exists positive constants $\kappa_H$ and $\beta>0$ such that for all $s,t\in T$
	$	|\sigma(t)-\sigma(s)|+|\mu(t)-\mu(s)|\leq \kappa_H|s-t|^{\beta}.$

 \end{itemize}

Denote by $\mathcal{C}(\mu_l,\mu_u, \sigma_l,\sigma_u, \beta,\kappa_H)$ the class of functions $\sigma(\cdot)$ and $\mu(\cdot)$ that satisfy Assumption C1.
We introduce the following uniform efficiency criterion.

\begin{definition}[Uniform polynomially efficient change of measure] We say that a change of measure $Q$  is uniformly polynomially efficient with the order $q\geq 0$ if there exists a constant $b_0\geq 0$ such that the importance sampling estimator
$$L_{\sigma,\mu}(b)=\frac{dP}{dQ}I{\Big(\sup_{t\in T} \{\sigma(t)f(t) +\mu(t)\}>b\Big)}$$
satisfies
\begin{equation}\label{logeff}
\sup_{b\geq b_0,\mu,\sigma\in \mathcal{C}(\mu_l,\mu_u, \sigma_l,\sigma_u, \beta,\kappa_H) } \frac{ Var (L_{\sigma,\mu}(b))}{ |\log w_{\sigma,\mu}(b) |^q w^{2}_{\sigma,\mu}(b)}<\infty.
\end{equation}
\end{definition}

Similar to the previous discussion, we consider the relative accuracy of a class of the importance sampling  estimators corresponding to a uniformly polynomially efficient change of measure. Let the $Q$ be uniformly polynomially efficient for $\sigma(\cdot),\mu(\cdot)\in \mathcal{C}(\mu_l,\mu_u, \sigma_l,\sigma_u, \beta,\kappa_H)$. Then, according to \eqref{eq:Cheby}, there exists some $\kappa_u>0$, such that the averaged estimator $Z_{\sigma,\mu}(b)=\frac{1}{n}\sum_{i=1}^n L^{(i)}_{\sigma,\mu}(b)$ based on $n=\kappa_u b^{2q}\delta^{-1}\varepsilon^{-2}$  i.i.d. Monte Carlo samples  satisfies
$$
\sup_{(\sigma,\mu)\in \mathcal{C}(\mu_l,\mu_u, \sigma_l,\sigma_u, \beta,\kappa_H)}P\left( |Z_{\sigma,\mu}(b)-w_{\sigma,\mu}(b)|>\varepsilon w_{\sigma,\mu}(b)\right)<\delta.
$$

\begin{remark}
Although the current paper focuses on rare-event simulation for the extremes of Gaussian random fields, the uniform efficiency criterion as well as the proposed method can be easily extended to  other Gaussian-related rare-event problems, such as the exponential integrals of Gaussian random fields \cite[e.g.,][]{LiuXu12,LiuXuTomacs},  where  the mean and variance functions are unspecified and we are interested in estimating a family of tail probabilities. Moreover, the proposed method can be extended to the estimation of non-Gaussian tail probabilities. For instance,
in statistical hypothesis testing with data generated independently from certain distribution with unknown parameters that are of interest, it often needs to evaluate the test power/error probabilities for a range of model parameters as the sample size increase; see \cite{li2016chernoff} for an example.
\end{remark}

\begin{remark}
In the literature, a similar uniform efficiency definition has been proposed in \cite{glasserman2008uniformly} to design an algorithm that is asymptotically efficient uniformly for a family of probability sets when estimating the tail probabilities of sums of light tailed random variables. Differently from this study, the random variable parameters are assumed to be known in their case.
\end{remark}


%
%
%
%
%

\section{Uniformly Efficient Estimation}\label{main}
\subsection{Discrete case}\label{sec:u}

We start with the case when $T$ contains finite points and propose a new change of measure which gives a uniformly efficient estimator. We assume $T:= \{t_{1},\cdots, t_M\}$.
We describe the new measure $Q$ in two ways.
First, we specify the sampling scheme of $f$ under $Q$ and then provide its Radon-Nikodym derivative with respect to $P$.
Under the measure $Q$, $f(t)$ is generated according to the following algorithm.

\begin{algorithm}[H]\label{alg:UGsim}
\caption{Simulating $f(\cdot)$ under $Q$}
\label{holderalg}
\KwIn{$T= \{t_{1},\cdots, t_M\}$, $\delta_b=ab^{-1}$ for some constant $a>0$.}
Simulate a random variable $\varsigma$  with respect to some positive continuous density function $g$ on $[\sigma_l,\sigma_u+\delta_b^2]$\;
Simulate a random variable $\nu$  with respect to some positive continuous density function $h$ on $[\mu_l,\mu_u+\delta_b]$\;
Simulate a random variable $\tau$ uniformly over $T =\{t_{1},\cdots, t_M\}$\;
 Given the realized $\varsigma$, $\nu$  and $\tau$, simulate $f(\tau)$ conditional on $\varsigma f(\tau)+\nu>b$\;
 Given $(\tau,f(\tau))$, simulate the Gaussian process $\{f(t): t\neq \tau, t\in T\}$ from the original conditional  distribution under $P$.

 \KwOut{$f(t)$ for $t\in T$}
\end{algorithm}

For the measure $Q$ defined above, it is not hard to verify that $P$ and $Q$ are mutually absolutely continuous with the Radon-Nikodym derivative being
$$\frac{dQ}{dP}=\int_{\mu_l}^{\mu_u+\delta_b}\int_{\sigma_l}^{\sigma_u+\delta^2_b}\frac{\sum_{i=1}^M I(\varsigma f(t_i)+\nu>b)}
{M P(\varsigma f(t_1)+\nu>b)} g(\varsigma)h(\nu)d\varsigma d\nu.$$
This gives the importance sampling estimator
\begin{eqnarray}\label{eq:proposed-est}
	L_{\sigma,\mu}(b) &=&
\left(\int_{\mu_l}^{\mu_u+\delta_b}\int_{\sigma_l}^{\sigma_u+\delta^2_b}\frac{\sum_{i=1}^M I(\varsigma f(t_i)+\nu >b)}
{M P(\varsigma f(t_1)+\nu>b)}
 g(\varsigma)h(\nu)d\varsigma d\nu\right)^{-1}\notag\\
 &&\times I{({\sup}_{i: t_i\in T} \sigma(t_i)f(t_i)+\mu(t_i)>b)}.
\end{eqnarray}
Note that under $Q$, if ${\max}_{t_i\in T} \sigma(t_i)f(t_i)+\mu(t_i)>b$, then
$\varsigma f(t_i)+\nu >b$ holds for all $i$, $\varsigma> \max_{t_i\in T}\sigma(t_i)$ and $\nu>\max_{t_i\in T}\mu(t_i)$. Therefore, the change of measure is well defined. 

We take a closer look at the proposed change of measure $Q$ by comparing it with the measure $Q^{\dagger}$ discussed in Section~\ref{sec:non-u}. We can see that steps 1 and 2 of Algorithm~\ref{alg:nonu} requires the knowledge of the mean and variance function $\mu$ and $\sigma$. When $\mu$ and $\sigma$ are unknown, running Algorithm~\ref{alg:nonu} with a misspecified $\mu'$ and $\sigma'$ may cause inefficiency. The proposed Algorithm~\ref{holderalg} avoids this inefficiency by introducing  prior probability density functions $g$ and $h$. Intuitively, the proposed algorithm explores each possible values of mean and variance of the random field at a random index (steps 1-3), and is a hybrid scheme for all $\sigma(\cdot)$ and $\mu(\cdot)$ that take values in the support of $g$ and $h$.
The next proposition states the uniform efficiency of the proposed change of measure.
\begin{proposition}\label{prop:variance-discrete}
Let 
$L_{\sigma,\mu}(b)$ be defined in \eqref{eq:proposed-est}, then 
 there exist constants $b_0$ and $\kappa_p$, independent of $\sigma(\cdot),\mu(\cdot)$ and $b$ and  for $b\geq b_0$,
$$
\frac{ E^Q(L^2_{\sigma,\mu}(b))}{M^2 b^6 w^2_{\sigma,\mu}(b)}\leq \kappa_p
$$
 for  all $\mu$ and $\sigma$ satisfying C1.
\end{proposition}
Note that $|\log(w_{\sigma,\mu}(b))|= O(b^2)$.
Therefore, the above proposition gives the uniformly polynomial efficiency of $Q$ with the order $q=3$ for the discrete case.
\begin{remark}
The parameter $\delta_b$ in Algorithm~\ref{holderalg} is introduced to control the second moment of the importance sampling estimator. 
Otherwise, consider the case of constant variance $\sigma\in[\sigma_l,\sigma_u]$ and zero mean $\mu=0$.
Then for $\sigma$ taking the value of $\sigma_u$, denote the corresponding estimator by $L_{\sigma_u,N}(b)$  and the second moment of $L_{\sigma_u,N}(b)$ is lower bounded by
\begin{eqnarray*}
 E^Q[L^2_{\sigma_u,N}(b)]
 &=&  E^Q\left[\left(\frac{dP}{dQ}\right)^2;~ \max_i\sigma_u f(t_i)>b~\right]
= E\left[ \frac{dP}{dQ} ; ~\max_i\sigma_u f(t_i)>b~\right]
 \\
 &=& {E\left[\left(\int_{\mu_l}^{\mu_u}\int_{\sigma_l}^{\sigma_u}\frac{\sum_{i=1}^M
I(\varsigma f(t_i) >b)}
{M P(\varsigma f(t_1)>b)}
 g(\varsigma)h(\nu)d\varsigma d\nu\right)^{-1}; ~ \max_i\sigma_u f(t_i)>b~\right]}\\
  &\geq&    P(\sigma_l f(0)>b) P(\max_i\sigma_u f(t_i)>b)\\
  &&\times
  {E\left[\left(\int_{\mu_l}^{\mu_u}\int_{\sigma_l}^{\sigma_u}
I(\max_i  f(t_i) >\varsigma^{-1}b)
 g(\varsigma)h(\nu)d\varsigma d\nu\right)^{-1}\Big| \max_i f(t_i)> \sigma_u^{-1}b~\right]}.
\end{eqnarray*}
However, the conditional expectation cannot be controlled and we have the estimator $L_{\sigma_u,N}(b)$ is not efficient for $\sigma=\sigma_u$.
\end{remark}

\begin{remark}
To  evaluate the Radon-Nikodym derivative in \eqref{eq:proposed-est},
we need to calculate the integral
$$\int_{\mu_l}^{\mu_u+\delta_b}\int_{\sigma_l}^{\sigma_u+\delta^2_b}\frac{\sum_{i=1}^M I(\varsigma f(t_i)+\nu >b)}
{M P(\varsigma f(t_1)+\nu>b)}
 g(\varsigma)h(\nu)d\varsigma d\nu.$$ Define
 \begin{equation}
 	l(z) = \int_{\mu_l}^{\mu_u+\delta_b}\int_{\sigma_l}^{\sigma_u+\delta^2_b}\frac{ I(\varsigma z+\nu >b)}
{\bar{\Phi}((b-\nu)/\varsigma)}
 g(\varsigma)h(\nu)d\varsigma d\nu,
 \end{equation}
 where $\bar{\Phi}(x)= \int_x^{\infty} \frac{1}{\sqrt{2\pi} }e^{-\frac{t^2}{2}}dt$ is the right tail probability of a standard Gaussian distribution,
 then we have
 \begin{equation*}
 	\int_{\mu_l}^{\mu_u+\delta_b}\int_{\sigma_l}^{\sigma_u+\delta^2_b}\frac{\sum_{i=1}^M I(\varsigma f(t_i)+\nu >b)}
{M P(\varsigma f(t_1)+\nu>b)}
 g(\varsigma)h(\nu)d\varsigma d\nu
 = \frac{1}{M}\sum_{i=1}^M l(f(t_i)).
 \end{equation*}
 Therefore, we only need to evaluate $l(f(t_i))$ for all $f(t_i)$ simulated by Algorithm~\ref{holderalg}. We use the following simplification for the function $l(z)$. Let $s = \frac{b-\nu}{\varsigma} $, then
 \begin{equation}\label{eq:lz-cv}
 \begin{split}
 		l(z) =& \int\int_{b -\varsigma s \in I_1,\varsigma\in I_2, s < z} \varsigma/ \bar{\Phi}(s) g(\varsigma) h(b-s\varsigma)d\varsigma ds
 		 = \int_{s < z } 1/\bar{\Phi}(s)\int_{\varsigma\in (\frac{b}{s}-\frac{1}{s}I_1)\cap I_2} \varsigma   h(b-s\varsigma)g(\varsigma)d\varsigma ds
 \end{split}
  \end{equation}
where   $I_1 = [\mu_l,\mu_u+\delta_b]$, and $I_2 = [\sigma_l,\sigma_u+\delta_b^2]$.
We can then choose $h(\cdot)$ and $g(\cdot)$ so that the inner integral in \eqref{eq:lz-cv} has a closed form expression. In particular, in the numerical examples  in this paper, we choose $g(\cdot)$ and $h(\cdot)$ to be the density functions of uniform distributions. In this case, let $
r(s) = \frac{1}{2}(\sigma_u+\delta_b^2-\sigma_l)^{-1}(\mu_u+\delta_b-\mu_l)^{-1} \int_{\varsigma\in (\frac{b}{s}-\frac{1}{s}I_1)\cap I_2} d\varsigma^2
$, then
$l(z)$ can be further simplified as
$$
l(z) = \int_{-\infty }^z r(s)/\bar{\Phi}(s)ds,
$$
which is a one-dimensional integral and can be evaluated numerically.
\end{remark}

\subsection{Continuous case}\label{Cont}
Direct simulation of a continuous random field is typically not a feasible task, and the change of measure proposed in the previous subsection is not directly applicable. 
Thus, we  use a discrete object to approximate the continuous fields for the implementation. 
In particular, we create a regular lattice covering $T$ in the following way. Let $G_{N,d}$ be a countable subset of $R^{d}$:
$G_{N,d} = \left\{\left(\frac{i_{1}}{N},\frac{i_{2}}{N},...,\frac{i_{d}}{N}\right): i_{1},...,i_{d}\in \mathbb Z\right\}.$
That is, $G_{N,d}$ is a regular lattice on $R^{d}$.
Furthermore, let
\begin{equation}\label{TM}
T_{N} = G_{N,d}\cap T,
\end{equation}
which is the sub-lattice intersecting with $T$.
Since $T$ is compact, $T_{N}$ is a finite set. 
We enumerate the elements in $T_{N} = \{t_{1},\cdots, t_M\}.$ Because $T$ is compact, we have $M=O(N^{d})$.
 Let $$w_{\sigma,\mu,N}(b) =P\left(\sup_{t_i\in T_N}\sigma(t_i)f(t_i)+\mu(t_i)>b\right). $$
 We use $w_{\sigma,\mu,N}(b)$ as a discrete approximation of $w_{\sigma,\mu}(b)$. We estimate $w_{\sigma,\mu,N}(b)$ by importance sampling, which is based on the change of measure proposed in  Section~\ref{sec:u}. In particular we define $Q_N$ and $P_N$ as the discrete versions (on $T_N$) of $Q$ and $P$  respectively. Then $dQ_N/dP_N$ takes the form:
\begin{equation}\label{disRN}
\frac{dQ_N}{dP_N}
=\int_{\mu_l}^{\mu_u+\delta_b}\int_{\sigma_l}^{\sigma_u+\delta^2_b}\frac{\sum_{i=1}^M
I(\varsigma f(t_i)+\nu>b)}
{M P(\varsigma f(t_1)+\nu>b)} g(\varsigma)h(\nu)d\varsigma d\nu.
\end{equation}
Note that here $M$ depends on $N$ and goes to infinity as $N\to\infty$.
This gives importance sampling estimator
\begin{eqnarray*}
	L_{\sigma,\mu,N}(b) &:=& \left(\int_{\mu_l}^{\mu_u+\delta_b}\int_{\sigma_l}^{\sigma_u+\delta^2_b}\frac{\sum_{i=1}^M
 I(\varsigma f(t_i)+\nu >b)}
{M P(\varsigma f(t_1)+\nu>b)}
 g(\varsigma)h(\nu)d\varsigma d\nu\right)^{-1}\\
 && \times I{({\sup}_{i: t_i\in T_N} \sigma(t_i)f(t_i)+\mu(t_i)>b)}.
\end{eqnarray*}
The discretization  usually introduces bias. The next two theorems control the bias and variance of the estimator $L_{\sigma,\mu,N}(b)$
under the following assumptions.
\begin{itemize}
	\item [C2] There exists a positive constant $\kappa_m$ such that 
$		\sup_{t\in T}\min_{t'\in T_N}|t-t'|\leq \frac{\kappa_m}{N}
$	for all $N$.
	\item [C3] The Gaussian random field $f$ is almost surely continuous.
	\item [C4] Define the correlation function 
 $
 r(s,t)=E(f(s)f(t)).
 $
 There exists $\beta'>0$ and $\kappa_H'>0$ such that
	\begin{equation}
		|r(t,s)-r(t',s')|\leq \kappa_H'[|t-t'|^{\beta'}+|s-s'|^{\beta'}]
	\end{equation}
	for all $s,t,s',t'\in T$. 
	\end{itemize}
\begin{theorem}\label{bias}
Let $\beta^*=\min(\beta,\beta')$ and $N_0(\varepsilon,b)=b^{2/\beta^*(\frac{3d}{\beta^*}+2-\varepsilon_0)}\varepsilon^{-2/\beta^*+\varepsilon_0}$.
Under Assumptions C1-C4, for any $\varepsilon_0>0$, there exist constants $\kappa_0$ and $b_0$ such that for any $\varepsilon\in (0,1)$, if 
$N\geq N_0(\varepsilon,b)$ and $b>b_0$, then
$$
\frac{\left| w_{\sigma,\mu,N}(b)- w_{\sigma,\mu}(b)\right|}
{w_{\sigma,\mu}(b)}< \varepsilon$$
uniformly for  $\mu,\sigma\in \mathcal{C}(\mu_l,\mu_u, \sigma_l,\sigma_u, \beta,\kappa_H)$.
\end{theorem}
\begin{theorem}\label{variance}
Let  $N_0(\varepsilon,b)$ be defined in  Theorem \ref{bias}.
Under Assumptions C1-C4, if $N\geq  N_0(\varepsilon,b)$, then there exist constants $b_0>0$ (depending on $\varepsilon_0$) and $\kappa_c>0$ such that
\begin{eqnarray*}
\sup_{b\geq b_0, \varepsilon\in (0,1)} \frac{E^{Q_N} L^2_{\sigma,\mu,N}(b)}
{b^q w_{\sigma,\mu}^{2}(b) \varepsilon^{-q_1}}<\kappa_c
\end{eqnarray*}
uniformly for  $\mu,\sigma\in \mathcal{C}(\mu_l,\mu_u, \sigma_l,\sigma_u, \beta,\kappa_H)$ with  $q=4d/\beta^*(\frac{3d}{\beta^*}+2+\varepsilon_0)+6$ and $q_1=4d/\beta^*+2d\varepsilon_0$.
\end{theorem}

We consider the relative accuracy of the importance sampling estimator based on $Q_N$. Let $L_{\sigma,\mu,N}^{(i)}(b)$ be i.i.d. copies of $L_{\sigma,\mu}(b)$ for $i=1,..,n$. Let 
\begin{equation}\label{eq:zn}
	Z_{\sigma,\mu,N}(b)=\frac{1}{n}\sum_{i=1}^n L_{\sigma,\mu,N}^{(i)}(b).
\end{equation}
With the aid of Chebyshev's inequality, we have
\begin{equation*}
	P(|Z_{\sigma,\mu,N}(b)-w_{\sigma,\mu}(b)|>\varepsilon w_{\sigma,\mu}(b) )\leq \frac{E(Z_{\sigma,\mu,N}(b)-w_{\sigma,\mu}(b))^2}{\varepsilon^2 w^2_{\sigma,\mu}(b)}.
\end{equation*}
The mean squared error $E(Z_{\sigma,\mu,N}(b)-w_{\sigma,\mu}(b))^2$ can be written as
\begin{equation*}
	E(Z_{\sigma,\mu,N}(b)-w_{\sigma,\mu}(b))^2 = [EZ_{\sigma,\mu,N} (b)- w_{\sigma,\mu}(b)]^2 + Var(Z_{\sigma,\mu,N}(b)).
\end{equation*}
The first and second terms on the right-hand side of the above display is the squared bias and the variance of the estimator $Z_{\sigma,\mu,N}(b)$, respectively.  If we choose $N=N_0(\varepsilon \delta^{1/2},b)$  according to Theorem~\ref{bias} and let $n=2\kappa_cb^{q}\varepsilon^{-q_1-2}\delta^{-\frac{q_1}{2}-1}$ where $q$ and $q_1$ are defined in Theorem~\ref{variance},
then the MSE is well controlled relative to $w_{\sigma,\mu}(b)$ and so is the relative accuracy. We summarize this result in the next corollary.


\begin{coro}\label{coro:accu}
	Under the Assumption C1-C4, let $Z_{\sigma,\mu,N}(b)$ be defined in \eqref{eq:zn}. If we choose $n= 2\kappa_c b^{q}\varepsilon^{-q_1-2}\delta^{-\frac{q_1}{2}-1}$ and $N=N_0(\varepsilon\delta^{1/2},b)
$, then
	\begin{equation}\label{eq:accuracy}
	 	P\left(|Z_{\sigma,\mu,N}(b)/w_{\sigma,\mu}(b)-1|>\varepsilon \right)<\delta.
	 \end{equation} 
\end{coro}
\begin{remark}
The  computational complexity for generating $Z_{\sigma,\mu,N}(b)$ is $n$ multiplied by the cost for generating one copy of $L_{\sigma,\mu,N}(b)$. The cost for generating $L_{\sigma,\mu,N}(b)$ is of order $O(M^3)=O(N^{3d})$, which is mainly the cost of generating a multivariate Gaussian vector (line 5 of Algorithm~\ref{holderalg}). The overall computational cost is also a polynomial in $\varepsilon,\delta$ and $b$. Algorithm with such a computation cost to achieve \eqref{eq:accuracy} is sometimes referred to as a fully polynomial randomized approximation scheme (FPRAS), see \cite{ABL09} for more details.
\end{remark}

\section{Simulation Studies}\label{Sim}
In this section, we present  numerical examples to show the performance of the proposed algorithm.
All the results are based on $n=10^4$ independent simulations.
The discretization size is chosen as $M=40$ in Example~\ref{example:cosine}-\ref{example:combined}. In each numerical example, we report the estimated tail probabilities, which will be referred to as ``est.'', along with the estimated standard deviations, that is  $sd^Q\{L_{\sigma,\mu}(b)\}=\sqrt{Var^Q\{L_{\sigma,\mu}(b)\}}$, which will be referred to as ``sd.''. The standard error of the estimator with $10^4$ Monte Carlo samples is sd.$/100$.
We also report the coefficient of variation (CV)  of the estimators, which is the ratio sd./est. of the estimators.

\bigskip
We start with the discrete setting in Example \ref{example:iid}, where $T=\{1,\cdots,100\}$ and
 $\{f(t), t=1,\cdots,100\}$ are i.i.d. standard normal random variables. We take $\mu(t)=0$ and $\sigma(t)=\sigma$ with $\sigma \in [0.3,1]$ for all $t\in T$, and the probability of interest is $P(\sigma\max_t f(t) >b)$ for   $b=3$. 
Table~\ref{result} displays the simulation results for
$\sigma = 0.3, 0.6$  and $1$ using the proposed method. 
 For different $\sigma$ values,  the estimates are close to the true values. Compared with the result of Algorithm 1 in Table \ref{result2}, the proposed method  gives   better  overall performance.

\begin{table}[h]
\centering
\begin{tabular}{ccccc}\hline
$\sigma$ &est.&sd.& CV  & Theoretical Value\\
\hline
 0.3 &  7.55e-22 & 5.33e-21 &  7.05 & 7.62e-22\\
   0.6&  2.93e-05 & 1.33e-04 & 4.52 & 2.87e-05\\
 1 &  1.26e-01 & 5.92e-01 & 4.69  & 1.26e-01\\
 \hline 
\end{tabular}
\caption{Estimates of $w_\sigma(b)$, $sd^Q(L_{\sigma,\mu}(b))$, and $sd^Q(L_{\sigma,\mu}(b))/w_\sigma(b)$. All results are based on $10^4$ independent simulations and thus the standard errors of the estimates are $sd^Q(L_{\sigma,\mu}(b))/100$.}
\label{result}
\end{table}
 
 We proceed to an example of a  continuous  Gaussian random field, whose tail probability of the supremum  is in a closed-form.
\begin{example}\label{example:cosine}
Consider the Gaussian random field $f(t)=X \cos t+Y\sin t,$ , where $X$ and $Y$ are independent standard Gaussian variables and $T=[0,3/4]$. We let $b=4$ and consider the class of constant variance and mean functions: $\sigma(t) = \sigma$ and $\mu(t) = \mu$, with $\sigma\in [0.5,1]$  and $\mu\in[-0.5,0.5]$. 
 \end{example}

For constant mean and variance functions considered in this example,
 the probability $P(\sup_{t\in T}(\sigma f(t)+\mu)>b)$ is known to be in a closed form \cite{AdlTay07}:
\begin{equation}\label{Truevalue}
  P\left(\sup_{0\leq t\leq 3/4 }(\sigma f(t)+\mu)>b\right)=\bar{\Phi}( (b-\mu)/\sigma )+\frac{3}{8\pi}e^{-(b-\mu )^2/(2\sigma^2)}.
\end{equation}
The simulation results for Example~\ref{example:cosine} are summarized in Table~\ref{table:sim1}. Similar to Example~\ref{example:iid}, we report the estimated probability, the standard deviation of the estimator, and its coefficient of variation. The theoretical value is computed according to \eqref{Truevalue}. We can see that for all combinations of $\sigma$ and $\mu$ in Table~\ref{table:sim1}, the estimated probabilities are close to the theoretical values. We also see that as the probability of interest decrease from $8.18 \times 10^{-6}$ to $4.01\times 10^{-12}$, the CV of the estimator does not increase substantially (from $2.7$ to $6.2$). This finding is consistent with our theoretical efficiency analysis of the proposed estimator.
\begin{table}[h!]
\centering
\begin{tabular}{cccccc}
\hline
$\sigma$ & $\mu$   & est.     & sd.       & CV  & theoretical value     \\ \hline
0.5   & 0.5  & 4.18E-12 & 2.59E-11 & 6.2 & 4.01E-12 \\
0.6   & 0.3  & 1.03E-09 & 4.38E-09 & 4.2 & 1.01E-09 \\
0.7   & 0.1  & 3.34E-08 & 1.18E-07 & 3.5 & 3.43E-08 \\
0.8   & -0.1 & 3.68E-07 & 1.19E-06 & 3.2 & 3.85E-07 \\
0.9   & -0.3 & 2.10E-06 & 5.97E-06 & 2.8 & 2.20E-06 \\
1     & -0.5 & 8.11E-06 & 2.20E-05 & 2.7 & 8.18E-06 \\ \hline
\end{tabular}
\caption{Simulation result for Example~\ref{example:cosine} with $b=4$ and $\delta_b=\frac{1}{b}$. Theoretical values are computed according to \eqref{Truevalue}.}
\label{table:sim1}
\end{table}

\begin{figure}[!ht]
\centering
\subfigure[Estimates as a function of $\beta_1$]{
      \includegraphics[scale=0.35]{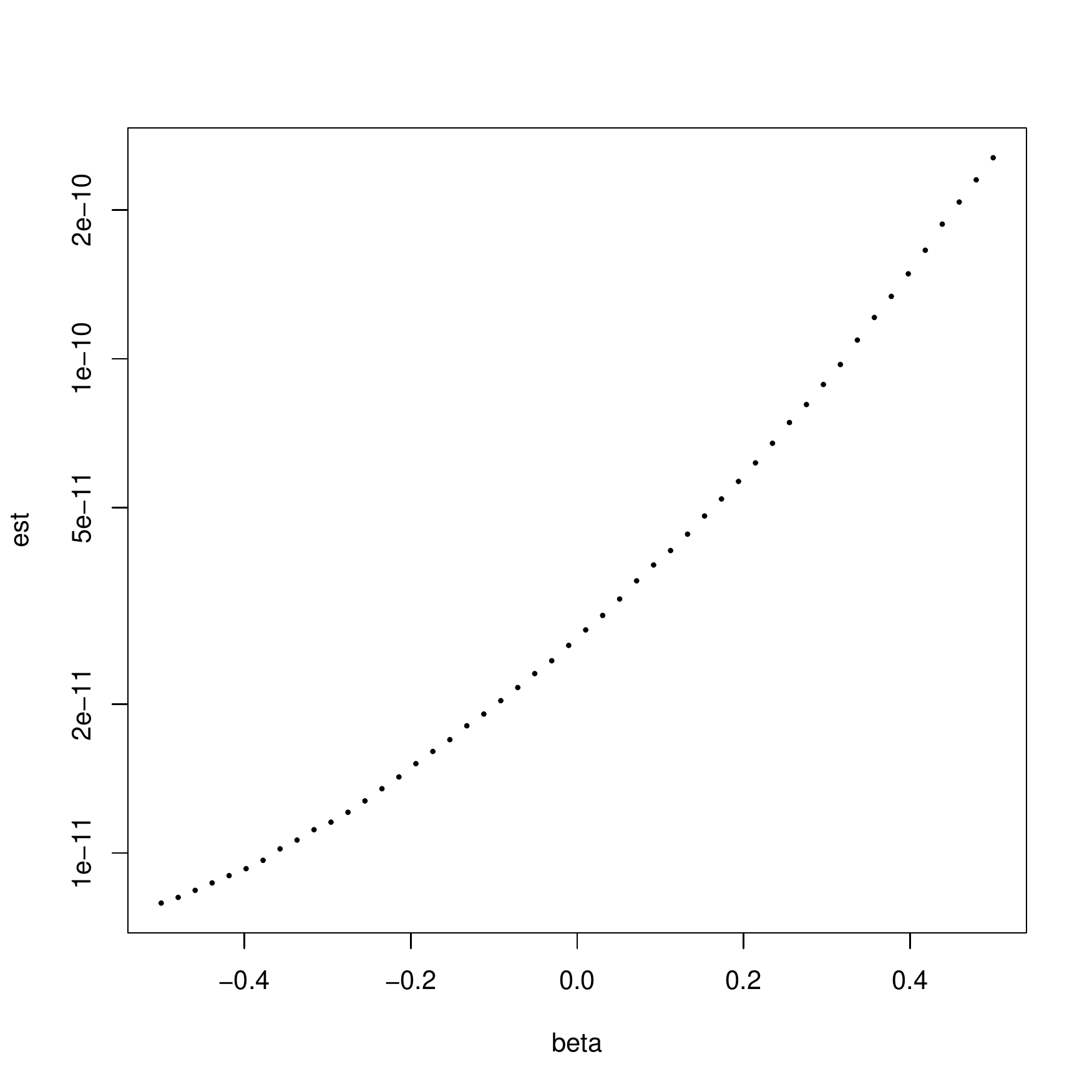}
    \label{fig:sim2-est}
}
\subfigure[CV as a function of $\beta_1$]{
      \includegraphics[scale=0.35]{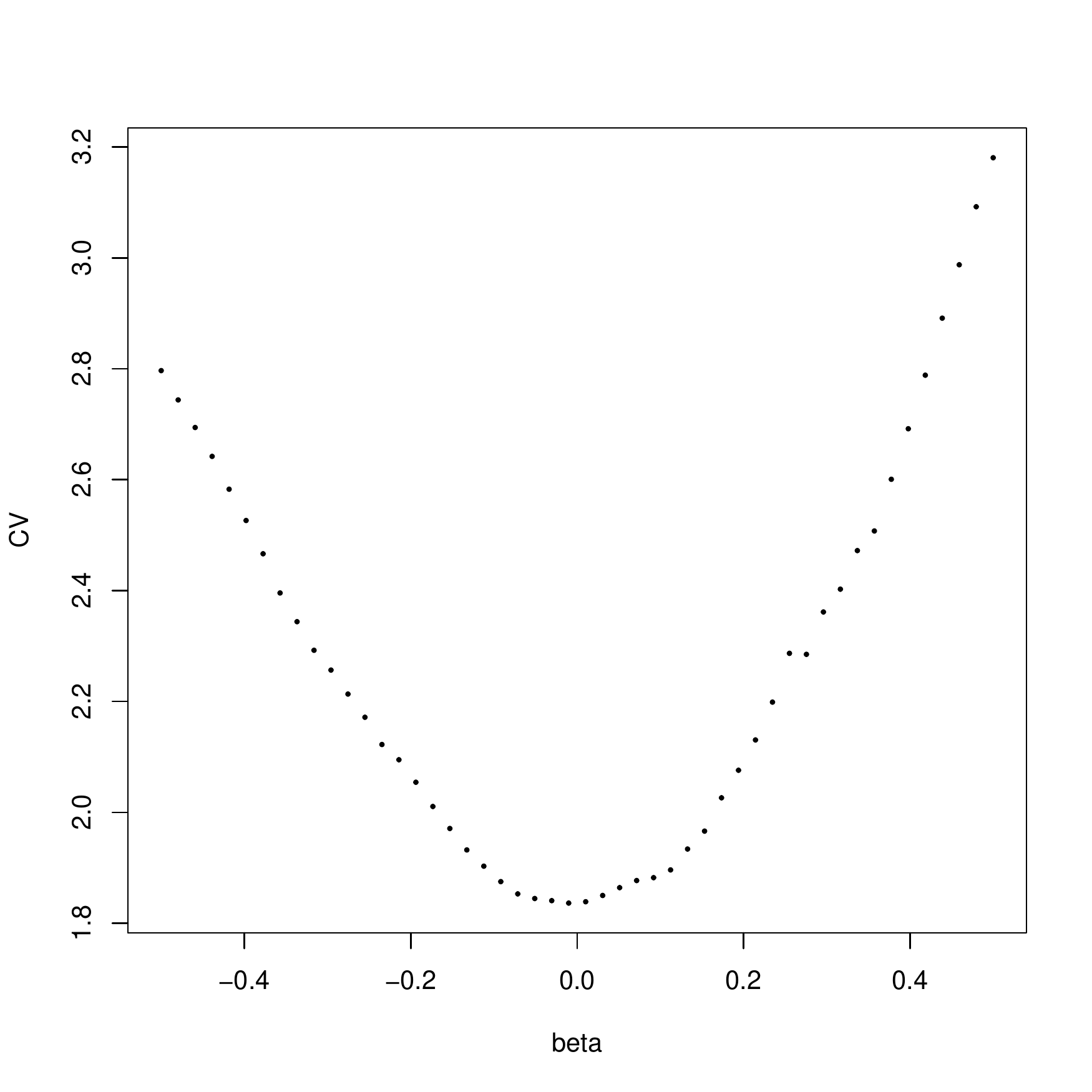}
    \label{fig:sim2-cv}
}
\caption[Optional caption for list of figures]{Simulation results for Example~\ref{example:linear}, where $b = 7$ and $\delta_b=1/b$.}
\label{fig:sim2}
\end{figure}

We proceed to examples where the mean and variance functions are not constants. We consider a continuous and centered Gaussian random field $\{f(t): 0\leq t \leq 1 \}$, whose covariance function is
\begin{equation}\label{eq:cov-abs}
	r(s,t)=E(f(s)f(t))=e^{-|s-t|}.
\end{equation}
In particular, in Example~\ref{example:linear} we consider a Gaussian random field with nonconstant mean and constant variance; in Example~\ref{example:variance} we consider a Gaussian field with constant mean and nonconstant variance; and in Example~\ref{example:combined},  both mean and variance functions are  nonconstant.
\begin{example}\label{example:linear}
Consider the Gaussian random field $f(t)$ defined in \eqref{eq:cov-abs}, and the class of variance and mean functions $\sigma(t) = 1$ and $\mu(t)=\beta_1 t$, for $\beta_1 \in [-0.5,0.5]$. The probability of interest is
$	P\left(
	\sup_{t\in[0,1]} f(t)+\beta_1 t>b
	\right)
$
for $b=7$.
\end{example}
 We summarize the simulation results for Example~\ref{example:linear} in Figure~\ref{fig:sim2}. Figure~\ref{fig:sim2-est} shows the scatter plot of the estimated probability (y-axis) against $\beta_1$ (x-axis). Figure~\ref{fig:sim2-cv} shows the scatter plot of the CV (y-axis)  against $\beta_1$ (x-axis). We see that the probability of interest is an increasing function in $\beta_1$. Moreover, when the estimated probability is within the range from $1\times 10^{-11}$ to $2\times 10^{-10}$, the CV of the estimator is always controlled within $3.2$, showing the good performance of the proposed estimation method.

\begin{figure}[ht]
\centering
\subfigure[Estimates as a function of $\beta_2$]{
      \includegraphics[scale=0.35]{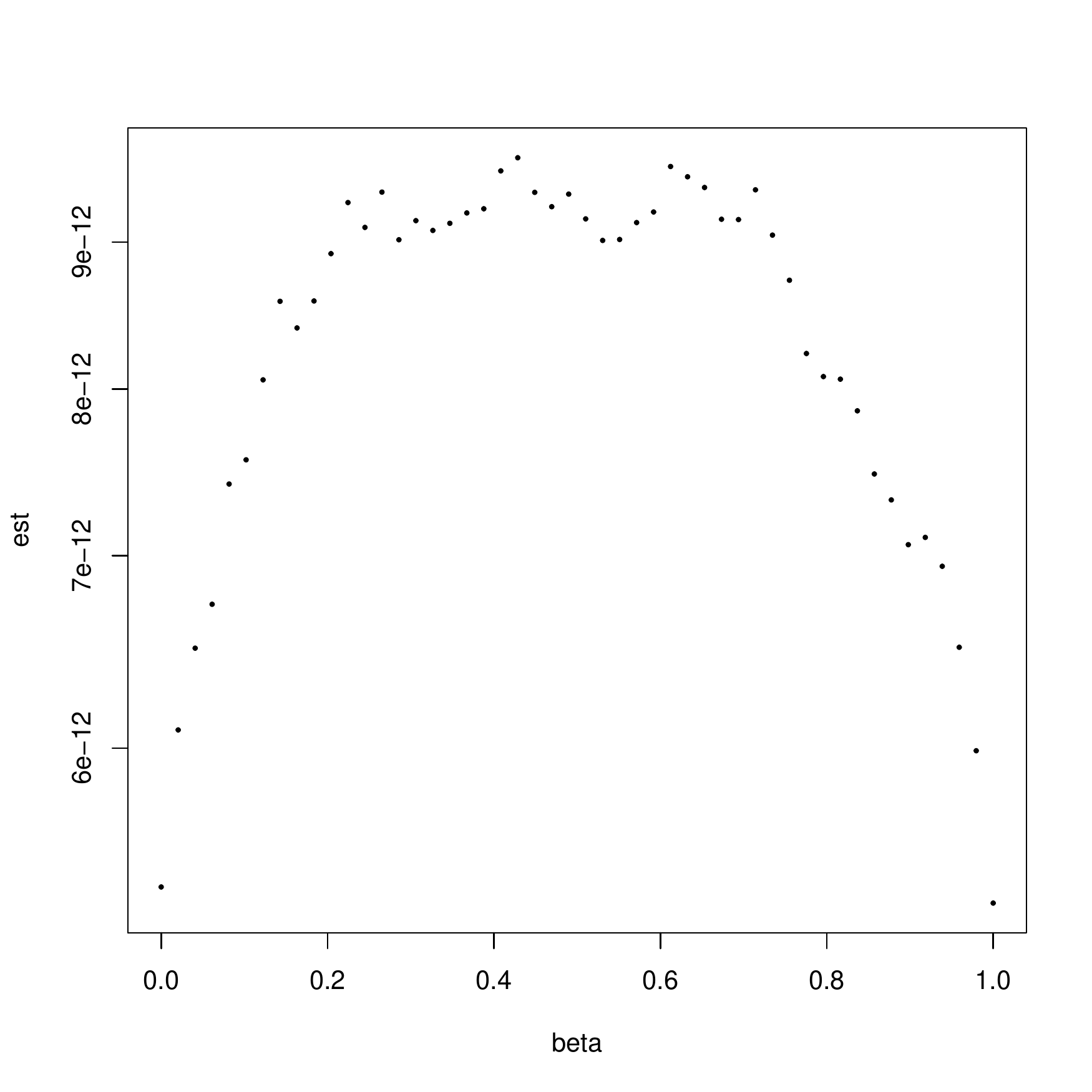}
    \label{fig:sim3-est}
}
\subfigure[CV as a function of $\beta_2$]{
      \includegraphics[scale=0.35]{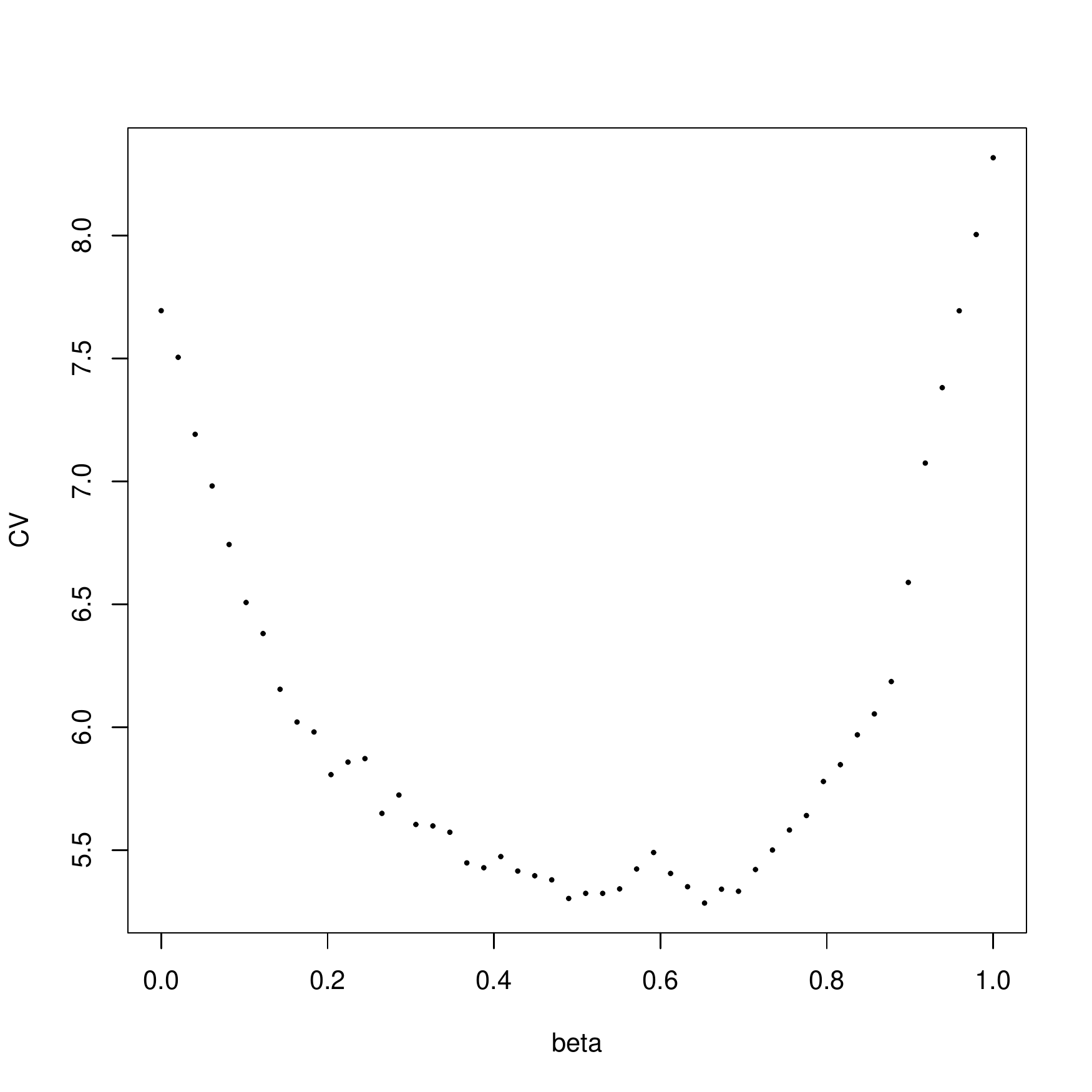}
    \label{fig:sim3-cv}
}
\caption[Optional caption for list of figures]{Simulation results for Example~\ref{example:variance}, where $b = 7$ and $\delta_b={1}/{b}$.}
\label{fig:sim3}
\end{figure}

\begin{example}\label{example:variance}
Consider the Gaussian random field $\{f(t), t\in T \}$ defined in \eqref{eq:cov-abs} and the class of variance and mean functions $\sigma(t) =  1- 0.5 (t-\beta_2)^2$ and $\mu(t)=0$, where $\beta_2 \in [0,1]$. The probability of interest is
$	P\left(
	\sup_{t\in[0,1]}  [ 1- 0.5 (t-\beta_2)^2] f(t)>b
	\right)
$ for $b=7$.
\end{example}
For Example~\ref{example:variance}, the scatter plot of estimated probability and the CV of the estimator are presented in Figure~\ref{fig:sim3}. Note that in Example~\ref{example:variance}, the maximum variance $\max_{t\in T}Var(\sigma(t)f(t))=Var(\sigma(\beta_2)f(\beta_2))=1$. Therefore, for all $\beta_2\in[0,1]$ the probability of interest has the same exponential decay rate $P(
	\sup_{t\in[0,1]}  \sigma(t) f(t)>b
	) = e^{-(1+o(1))\frac{b^2}{2\max_{t\in T}Var(\sigma(t)f(t))}}= e^{-(1+o(1))b^2/2}$, as $b\to\infty$. In Figure~\ref{fig:sim3-est}, we see that the estimated probability is relatively small when $\beta_2$ is close to the boundary values $0$ or $1$, compared to the case when $\beta_2\in[0.2,0.8]$ and is far away from the boundary values. For $\beta_2\in[0.2,0.8]$ the estimated probability stays around $9\times 10^{-12}$ and does not fluctuate much. For all $\beta_2\in[0,1]$, the maximum CV of the estimator is controlled within $10$. This is again consistent with our theoretical results.

\begin{example}\label{example:combined}
Consider the Gaussian random field $\{f(t),t\in T\}$ defined in \eqref{eq:cov-abs}, and the class of variance and mean functions $\sigma(t) =  1- 0.5 (x-\beta_2)^2$ and $\mu(t)=\beta_1 t$, where $\beta_1\in[-0.5,0.5]$ and $\beta_2 \in [0,1]$. The probability of interest is
$P(
	\sup_{t\in[0,1]} \{[1-0.5\times (t-\beta_2)^2]f(t) + \beta_1 t\}>b
	),$ for $b=7$.
\end{example}
Table~\ref{table:sim4} shows the simulated results for different choices of $\beta_1$ and $\beta_2$. We see that the estimated probabilities range from $4.2\times 10^{-12}$ to $1.16\times 10^{-10}$. The maximum CV in Table~\ref{table:sim4} is $9.9$. This means that the standard error of the averaged Monte Carlo estimator with $10^4$ samples is controlled within $9.9\%\times E^Q L_{\sigma,\mu}(b)$. 
\begin{table}[ht]
\centering
\begin{tabular}{ccccc}
\hline
$\beta_1$ & $\beta_2$ & est.      & sd.       & CV  \\ \hline
-0.50 & 0.00  & 4.20E-12 & 4.03E-11 & 9.6 \\
-0.33 & 0.17  & 5.60E-12 & 3.69E-11 & 6.6 \\
-0.17 & 0.33  & 5.69E-12 & 3.29E-11 & 5.8 \\
0.00  & 0.50  & 8.78E-12 & 5.09E-11 & 5.8 \\
0.17  & 0.67  & 2.09E-11 & 1.27E-10 & 6.1 \\
0.33  & 0.83  & 5.82E-11 & 4.04E-10 & 6.9 \\
0.50  & 1.00  & 1.16E-10 & 1.15E-09 & 9.9 \\ \hline
\end{tabular}
\caption{Simulation results for Example~\ref{example:combined}, where $b=7$ and $\delta_b=2/b$.}
\label{table:sim4}
\end{table}
\section{Proofs of main results}\label{SecProof}
Throughout the proof, we write $a(b)=O(c(b))$ if there exists a positive constant $\kappa$, independent of $b,\sigma(\cdot),\mu(\cdot)$, such that $|a(b)|/|c(b)|\leq \kappa$. We also write $a(b)=o(c(b))$ if $|a(b)|/|c(b)|\to 0$ as $b\to\infty$, uniformly in $\sigma(\cdot)$ and $\mu(\cdot)$ satisfying Assumption C1.
We will use $\tkappa$ as a generic notation to denote large and not-so-important
constants (independent of $\mu,\sigma$ and $b$) whose value may vary from place to place. Similarly, we use $\tep$ as a generic notation for
small positive constants.

\begin{proof}[\bf Proof of Proposition~\ref{prop:nonu}]
We start with the proof of Proposition~\ref{prop:nonu} (i). We can see that if ${\max}_{t_i\in T}\sigma'(t_i)f(t_i)+\mu(t_i)>b$, then $\max_{t_i\in T}\sigma(t_i) f(t_i)+\mu(t_i)>b$ always happens and the change of measure is well defined.
We have
\begin{eqnarray*}
&&E^{Q^{\dagger}}\left[\left(\frac{dP}{dQ^{\dagger}}\right)^{2};
{\max}_{t_i\in T}\sigma'(t_i)f(t_i)+\mu(t_i)>b\right]\\
&=& E\left[\frac{dQ^{\dagger}}{dP}\times\left(\frac{dP}{dQ^{\dagger}}\right)^{2};
{\max}_{t_i\in T}\sigma'(t_i)f(t_i)+\mu(t_i)>b\right]\\
&=& E\left[\frac{dP}{dQ^{\dagger}};
{\max}_{t_i\in T}\sigma'(t_i)f(t_i)+\mu(t_i)>b\right]\\
&=&   {E\left[\frac{\sum_{i=1}^M P(\sigma(t_i) f(t_i)+\mu(t_i)>b)}
{\sum_{i=1}^M I{(\sigma(t_i) f(t_i)+\mu(t_i)>b)}};
{\max}_{t_i\in T}\sigma'(t_i)f(t_i)+\mu(t_i)>b\right]}.
\end{eqnarray*}
Because $\sum_{i=1}^M I{(\sigma(t_i) f(t_i)+\mu(t_i)>b)}\leq M$, the above display is further bounded from below by
\begin{eqnarray*}
&=&
\frac{1}{M}
\left(\sum_{i=1}^M P(\sigma(t_i) f(t_i)+\mu(t_i)>b)\right)\times w_{\sigma',\mu}(b)\\
&\geq&
\frac{1}{M}\max_{t_i\in T} P(\sigma(t_i) f(t_i)+\mu(t_i)>b)\times w_{\sigma',\mu}(b)\\
&=& \exp\left\{-(1+o(1))\frac{b^2}{2\max_{t_i\in T} \sigma(t_i)^2}-(1+o(1))\frac{b^2}{2\max_{t_i\in T} \sigma'(t_i)^2}\right\},
\end{eqnarray*}
where we used the following lemma, whose proof is given  in Section 5.1,  to obtain that
$$w_{\sigma',\mu}(b)=\exp\left\{-(1+o(1))\frac{b^2}{2\max_{t_i\in T} \sigma'(t_i)^2}\right\}.
$$
\begin{lemma}\label{lemma:tail-approx}
	Let $\{f(t):t\in T \}$ be a centered, unit variance and continuous Gaussian random field living on a compact set $T$. Assume that $\sigma(t)>0$ and $\mu(t)$ are continuous functions. Then, there exists positive $\tep$ such that 
	$$
	P\Big(\sup_{t\in T}\sigma(t)f(t)+\mu(t)>b\Big) = e^{-(1+o(1))\frac{b^2}{2\max_{t\in T}\sigma^2(t)}} \mbox{ and } P\Big(\sup_{t\in T}\sigma(t)f(t)+\mu(t)>b\Big) \geq \tep b^{-1}\max_{t\in T}e^{-\frac{(b-\mu(t))^2}{2\sigma^2(t)}}.
	$$
\end{lemma}
Under the assumption that
$\max_{t_i\in T}\sigma'(t_i)<\max_{t_i\in T}\sigma(t_i)$, we know that for $\varepsilon<\frac{1}{2}(1-\frac{\max \sigma'(t_i)}{\max \sigma(t_i)})$
$$
\frac{E^{Q^{\dagger}}\Big[\left(\frac{dP}{dQ^{\dagger}}\right)^{2};
{\max}_{t_i\in T}\sigma'(t_i)f(t_i)+\mu(t_i)>b\Big]}
{ w_{\sigma',\mu}^{2-\epsilon}(b)} 
\geq w_{\sigma',\mu}^{-\epsilon}(b),
$$
which tends to infinity as $b\to\infty$.

We proceed to the proof of part (ii). Let $t'_{\max}=\arg\max_{t\in T}\sigma'(t)$. We consider the event $F=\{
b/\sigma'(t'_{\max})<f(t'_{\max})<\min_{t_i\in T}[b/\sigma(t_i)]
\}$. Because $\max_{t_i\in T}\sigma'(t_i)>\max_{t_i\in T}\sigma(t_i)$, $F$ is non-empty and $F\subset \{
\max_{t_i\in T}\sigma'(t_i)f(t_i)+\mu'(t_i)>b
\}$. Moreover, according to the sampling scheme in Algorithm~\ref{alg:nonu}, we have $Q^{\dagger}(F)>0$. On the other hand, when the event $F$ happens, $\sum_{i=1}^M I{(\sigma(t_i) f(t_i)}>b)=0$, therefore $Q^{\dagger}(\frac{dP}{dQ^{\dagger}}=\infty)\geq Q^{\dagger}(F)>0$. In other word, $\frac{dP}{dQ^{\dagger}}$ is not well-defined.
\end{proof}
\bigskip
\begin{proof}[\bf Proof of Proposition~\ref{prop:variance-discrete}]
Define the random index $t^*:=\arg\max_{t\in T} [\sigma(t)f(t)+\mu(t)]$. We restrict our analysis to the integral over the region $[\mu(t^*),\mu(t^*)+\delta_b]\times [\sigma(t^*),\sigma(t^*)+\delta_b^2]$ and arrive at
\begin{align}
& E^Q[L_{\sigma,\mu}^2(b)]\notag\\
=&   {E^Q\biggr[\biggr(\int_{\mu_l}^{\mu_u+\delta_b}\int_{\sigma_l}^{\sigma_u+\delta_b^2}\frac{\sum_{i=1}^M
I(\varsigma f(t_i)+\nu >b)}
{M P(\varsigma f(t_1)+\nu>b)}
 g(\varsigma)d\varsigma d\nu\biggr)^{-2}; {\max}_{t_i\in T}~ \sigma(t_i)f(t_i)+\mu(t_i)>b\biggr]}\notag\\
\leq&
{E^Q\biggr[\biggr(\int_{\mu(t^*)}^{\mu(t^*)+\delta_b}\int_{\sigma(t^*)}^{\sigma(t^*)+\delta^2_b}
\frac{\sum_{i=1}^M
I(\varsigma f(t_i)+\nu >b) }{M P(\varsigma f(t_1)+\nu>b)}g(\varsigma)h(\nu)d\varsigma d\nu\biggr)^{-2}; {\max}_{t_i\in T}~\sigma(t_i)f(t_i)+\mu(t_i)>b\biggr]}\notag\\
=& {E^Q\biggr[\biggr(\int_{\mu(t^*)}^{\mu(t^*)+\delta_b}\int_{\sigma(t^*)}^{\sigma(t^*)+\delta^2_b}
\frac{\sum_{i=1}^M
I(\varsigma f(t_i)+\nu >b) }{M \bar{\Phi}(\frac{b-\nu}{\varsigma})}g(\varsigma)h(\nu)d\varsigma d\nu\biggr)^{-2}; {\max}_{t_i\in T}~\sigma(t_i)f(t_i)+\mu(t_i)>b\biggr]}\notag\\
\label{eq:bound-cont}
\end{align}
Note that for all $(\varsigma,\nu)\in [\mu(t^*),\mu(t^*)+\delta_b]\times [\sigma(t^*),\sigma(t^*)+\delta_b^2]$, we have $\varsigma f(t^*)+\nu\geq \max_{t_i\in T}\sigma(t_i)f(t_i)+\mu(t_i)$. Therefore, the event $ \max_{t_i\in T}\sigma(t_i)f(t_i)+\mu(t_i)>b$ implies $\varsigma f(t^*)+\nu\geq b$. Consequently, $\sum_{i=1}^M
I(\varsigma f(t_i)+\nu >b) \geq 1$ on the event ${\max_{t_i\in T}}\sigma(t_i)f(t_i)+\mu(t_i)>b$. Therefore, 
\eqref{eq:bound-cont} is further bounded from above by
\begin{eqnarray}
&\leq&M^2
{E^Q\biggr[\biggr(\int_{\mu(t^*)}^{\mu(t^*)+\delta_b}\int_{\sigma(t^*)}^{\sigma(t^*)+\delta^2_b}
\frac{g(\varsigma)h(\nu) }{ \bar{\Phi}(\frac{b-\nu}{\varsigma})}d\varsigma d\nu\biggr)^{-2}; {\max}_{t_i\in T}\sigma(t_i)f(t_i)+\mu(t_i)>b\biggr]}\notag\\
&\leq&O(1)M^2 E^Q\biggr[\Big(\int_{\mu(t^*)}^{\mu(t^*)+\delta_b}\int_{\sigma(t^*)}^{\sigma(t^*)+\delta^2_b}
{g(\varsigma)h(\nu)}{b e^{\frac{(b-\nu)^2}{2\varsigma^2}}
}d\varsigma d\nu\Big)^{-2}; {\max}_{t_i\in T} \sigma(t_i)f(t_i)+\mu(t_i)>b\biggr]\notag\\
\label{eq:cont2}
\end{eqnarray}
Note that for all $(\varsigma,\nu)\in[\mu(t^*),\mu(t^*)+\delta_b]\times [\sigma(t^*),\sigma(t^*)+\delta_b^2]$, we have
$ b e^{\frac{(b-\nu)^2}{2\varsigma^2}} = O(1)b e^{\frac{(b-\mu(t^*))^2}{2\sigma^2(t^*)}}$. Therefore, \eqref{eq:cont2} is bounded from above by
\begin{eqnarray}	
&\leq&O(1)M^2 E^Q\biggr[\Big(\int_{\mu(t^*)}^{\mu(t^*)+\delta_b}\int_{\sigma(t^*)}^{\sigma(t^*)+\delta^2_b}
{g(\varsigma)h(\nu)}{b e^{\frac{(b-\mu(t^*))^2}{2\sigma^2(t^*)}}
}d\varsigma d\nu\Big)^{-2}; {\max}_{t_i\in T} \sigma(t_i)f(t_i)+\mu(t_i)>b\biggr]\notag\\
&=& O(1)M^2 \delta_b^{-6} b^{-2} E^Q\Big[e^{-\frac{(b-\mu(t^*))^2}{ \sigma^2(t^*)}};{\max}_{t_i\in T} \sigma(t_i)f(t_i)+\mu(t_i)>b\Big]\notag\\
&\leq & O(1)M^2 \delta_b^{-6} b^{-2} \max_{t_i\in T}e^{-\frac{(b-\mu(t_i))^2}{ \sigma^2(t_i)}}.\label{eq:variance-control}
\end{eqnarray}
On the other hand, according to Lemma~\ref{lemma:tail-approx} we have
\begin{equation*}
	P\Big(\sup_{t_i\in T} \sigma(t_i) f(t_{i})+\mu(t_{t_i})>b\Big)\geq \tep b^{-1}\max_{t_i\in T}e^{-\frac{b-\mu(t_i)}{2\sigma^2(t_i)}}.
\end{equation*}
Combining this and \eqref{eq:variance-control}, we have that there exists $b_0$ sufficiently large such that for $b\geq b_0$
$$
\frac{ E^Q[L^2_{\sigma,\mu}(b); {\max}_{t_i\in T} \sigma(t_i)f(t_i)+\mu(t_i)>b]}{M^2 b^6 w_{\sigma,\mu}^2(b)}= O(1).
$$
This completes our proof.
\end{proof}

\medskip
\begin{proof}[\bf Proof of Theorem \ref{bias}]
Note that $\sup_{t\in T}\sigma(t)f(t)+\mu(t)\geq \sup_{t\in T_N}\sigma(t)f(t)+\mu(t)$, we have
\begin{eqnarray*}
	&&\Big|P\Big(\sup_{t\in T}\sigma(t)f(t)+\mu(t)>b\Big)- P\Big(\sup_{t\in T_N}\sigma(t)f(t)+\mu(t)>b\Big)\Big|\\
	&=& P\Big(\sup_{t\in T}\sigma(t)f(t)+\mu(t)>b, \sup_{t\in T_N}\sigma(t)f(t)+\mu(t)\leq b\Big).
\end{eqnarray*}
We  split the above probability into two parts.
\begin{eqnarray*}
	&&P\Big(\sup_{t\in T}\sigma(t)f(t)+\mu(t)>b, \sup_{t\in T_N}\sigma(t)f(t)+\mu(t)\leq b\Big)\\
	&=&
	P\Big(b<\sup_{t\in T}\sigma(t)f(t)+\mu(t)\leq b+\frac{\gamma}{b}, \sup_{t\in T_N}\sigma(t)f(t)+\mu(t)\leq b\Big)\\
	&&+
	P\Big(\sup_{t\in T}\sigma(t)f(t)+\mu(t)> b+\frac{\gamma}{b}, \sup_{t\in T_N}\sigma(t)f(t)+\mu(t)\leq b\Big),
\end{eqnarray*}
which is further bounded from above by
\begin{equation}\label{eq:two-terms}
	P\Big(b<\sup_{t\in T}\sigma(t)f(t)+\mu(t)\leq b+\frac{\gamma}{b}\Big)+ P\Big(\sup_{t\in T}\sigma(t)f(t)+\mu(t)> b+\frac{\gamma}{b}, \sup_{t\in T_N}\sigma(t)f(t)+\mu(t)\leq b\Big),
\end{equation}
where we will choose $\gamma$ later.
We proceed to upper bounds of the above two terms separately.
For the first term, we apply the following Lemma.
\begin{lemma}[Proposition 6.5 of \cite{ABL09}]\label{lemma:adler}
	Under Assumptions C1,C3 and C4, for any $v>0$, let $\beta^*=\min(\beta,\beta')$ and $\rho=\frac{2d}{\beta^*}+dv+1$, where $d$ is the dimension of $T$. There exists constants $b_0,\lambda\in (0,\infty)$ so that for all $b\geq b_0\geq 1$,
	\begin{equation}
		P\Big(
		\max_{t\in T} \sigma(t)f(t)+\mu(t)\leq b+\frac{\gamma}{b} ~\Big |~ \max_{t\in T}\sigma(t)f(t)+\mu(t)>b
		\Big)\leq \lambda a b^{\rho}.
	\end{equation}
\end{lemma}
With the aid of the above lemma with $v=\frac{1}{{\beta^*}}$, we have for $b\geq b_0$
\begin{eqnarray*}
	&&P\Big(b<\sup_{t\in T}\sigma(t)f(t)+\mu(t)\leq b+\frac{\gamma}{b}\Big)\\
	&=& P\Big(\max_{t\in T}\sigma(t)f(t)+\mu(t)>b\Big) P\Big(
		\max_{t\in T} \sigma(t)f(t)+\mu(t)\leq b+\frac{\gamma}{b} ~\Big |~ \max_{t\in T}\sigma(t)f(t)+\mu(t)>b
		\Big)\\
	&\leq &\lambda \gamma b^{\rho} P\Big(\max_{t\in T}\sigma(t)f(t)+\mu(t)>b\Big)
\end{eqnarray*}
with $\rho=\frac{3d}{{\beta^*}}+1$.
We choose $\gamma:=2^{-1}\lambda^{-1}b^{-\rho}\varepsilon$, then the above display gives the following upper bound for the first term in \eqref{eq:two-terms}
\begin{eqnarray*}
	P\Big(b<\sup_{t\in T}\sigma(t)f(t)+\mu(t)\leq b+\frac{\gamma}{b}\Big)\leq \frac{\varepsilon}{2} w_{\sigma,\mu}(b).
\end{eqnarray*}
We proceed to the second term in \eqref{eq:two-terms}. According to Assumption C2, we have
\begin{eqnarray*}
		&&P\Big(\sup_{t\in T}\sigma(t)f(t)+\mu(t)> b+\frac{\gamma}{b}, \sup_{t\in T_N}\sigma(t)f(t)+\mu(t)\leq b\Big)\\
&\leq& P\Big(
\sup_{t,s\in T, |t-s|\leq \kappa_m/{N}}|\sigma(t)f(t)+\mu(t)-(\sigma(s)f(s)+\mu(s))|>\frac{\gamma}{b}
\Big),
\end{eqnarray*}
which is further bounded from above by
\begin{equation}\label{eq:bound1}
\begin{split}
P\Big(
\sup_{t,s\in T, |t-s|\leq \kappa_m/{N}}|\sigma(t)f(t)-\sigma(s)f(s)|+\sup_{t,s\in T, |t-s|\leq \kappa_m/{N}}|\mu(t)-\mu(s)|>\frac{\gamma}{b}
\Big).
\end{split}
\end{equation}
According to Assumption C1, we have
\begin{equation*}
	\sup_{t,s\in T, |t-s|\leq \kappa_m/{N}}|\mu(t)-\mu(s)|= O(\kappa_m^{{\beta^*}}/{N^{{\beta^*}}}).
\end{equation*}
Plugging this into \eqref{eq:bound1}, we have
\begin{equation}
\begin{split}
	&P\Big(
\sup_{t,s\in T, |t-s|\leq \kappa_m/{N}}|\sigma(t)f(t)+\mu(t)-\sigma(s)f(s)+\mu(s)|>\frac{\gamma}{b}
\Big)\\
\leq & P\Big(
\sup_{t,s\in T, |t-s|\leq \kappa_m/{N}}|\sigma(t)f(t)-\sigma(s)f(s)|>\frac{\gamma}{b}- \kappa_m^{{\beta^*}}/{N^{{\beta^*}}}
\Big).
\end{split}
\end{equation}
We choose $N\geq \tkappa \lambda^{1/{\beta^*}} b^{(\rho+1)/{\beta^*}}\varepsilon^{-1/{\beta^*}}$ for $\tkappa$ sufficiently large, then $\frac{\gamma}{b}- \kappa_m^{{\beta^*}} \frac{1}{N^{{\beta^*}}}>\frac{\gamma}{2b}$. Therefore, we further have
\begin{equation}
\begin{split}
	&P\Big(
\sup_{t,s\in T, |t-s|\leq \kappa_m/{N}}|\sigma(t)f(t)+\mu(t)-\sigma(s)f(s)+\mu(s)|>\frac{\gamma}{b}
\Big)\\
\leq & P\Big(
\sup_{t,s\in T, |t-s|\leq \kappa_m/{N}}|\sigma(t)f(t)-\sigma(s)f(s)|>\frac{\gamma}{2b}
\Big).
\end{split}
\end{equation}
To control the above probability, we use the following lemma known as the Borell-TIS lemma, which is proved
independently by \cite{Bor75} and \cite{CIS}.
\medskip
\begin{lemma}[Borell-TIS]\label{LemBorell}
Let $\{f(t);t\in \mathcal U\}$, where $\mathcal U$ is a compact set,
be a mean zero Gaussian random field. $f$ is almost surely
bounded on $\mathcal U$. Then, $E[\sup_{\mathcal U}f(t)]<\infty,$ and
$
P\left(\sup_{t\in \mathcal{U}}f\left(  t\right)
-E[\sup_{t\in\mathcal{U}}f\left( t\right)  ]\geq b\right)\leq
\exp\left(
-\frac{b^{2}}{2\sigma_{\mathcal{U}}^{2}}\right) ,
$
where $\sigma_{\mathcal{U}}^{2}=\sup_{t\in \mathcal{U}}\hbox{Var}[f( t)].$
\end{lemma}

\medskip
We define a new Gaussian random field
\begin{equation}
	\xi(s,t)=\sigma(s)f(s) -\sigma(t)f(t).
\end{equation}
The next lemma, whose proof will be provided in Section~\ref{sec:proof-lemmas}, characterizes $E {\sup_{t,s\in T, |t-s|\leq \kappa_m/{N}}  } \xi(s,t)$.
\begin{lemma}\label{lemma:esup}
For all $\sigma,\mu$ and $f$ satisfying Assumptions C1, C3 and C4, there is a  uniform constant $\kappa_{\xi}>0$ such that
$$
E {\sup_{t,s\in T, |t-s|\leq \kappa_m/{N}}} |\xi(s,t)|<\kappa_{\xi} N^{-\beta^*/2} \log N
$$
\end{lemma}
Furthermore, the variance of $\xi(s,t)$ is bounded from above by
\begin{equation}\label{eq:var-xi}
	Var(\xi(s,t))= (\sigma(s)-\sigma(t))^2 +2\sigma(s)\sigma(t)(1-r(s,t))\leq \kappa_H^2|s-t|^{2{\beta^*}}+2\sigma_u^2 |s-t|^{{\beta^*}}\leq O( |s-t|^{{\beta^*}}).
\end{equation}
According to Assumption C1 and C4, the above display is further bounded from above by
\begin{equation}
	Var(\xi(s,t))\leq O( N^{-{{\beta^*}}})
\end{equation}
We choose $N$ such that $\kappa_{\xi} N^{-\frac{\beta^*}{2}}\log N \leq \frac{\gamma}{4b}$. Then according to the Borell-TIS lemma and Lemma~\ref{lemma:esup}, we have
\begin{equation}
	P\Big(
	\sup_{|t-s|\leq\frac{\kappa_{m}}{N}}|\xi(s,t)|>\frac{\gamma}{4b}
	\Big)\leq  \exp\left(
- \tep \frac{\gamma^2}{N^{-{\beta^*} }b^2}
\right).
\end{equation}
The above display is of order $o(\varepsilon w_{\sigma,\mu}(b))$ if
$
\frac{\gamma^2}{N^{-\beta^*}b^2}\geq \tkappa^{\beta^*} \max(-\log \varepsilon, b^2)
$, for a large enough and possibly different constant $\tkappa$. Therefore, it is sufficient to choose
$N\geq \tkappa \max(-\log \varepsilon, b^2)^{1/\beta^*}b^{2/\beta^*}\gamma^{-2/\beta^*}(\log b)^{\tkappa}$.
Combining this with our choice of $\gamma$, and recall our choice of $\rho$ in Lemma~\ref{lemma:adler} it is sufficient to choose
$N\geq \tkappa \max(-\log \varepsilon, b^2)^{1/\beta^*} b^{2/\beta^*+2/\beta^*(\frac{3d}{\beta^*}+1)}\varepsilon^{-2/\beta^*}(\log b)^{\tkappa}$, which is bounded by  $N_0=b^{2/\beta^*(\frac{3d}{\beta^*}+2+\varepsilon_0)}\varepsilon^{-2/\beta^*-\varepsilon_0}$ for any $\varepsilon_0>0$ and $b$ sufficiently large. This completes our proof.
\end{proof}

\medskip\begin{proof}[\bf Proof of Theorem~\ref{variance}]
According to Proposition~\ref{prop:variance-discrete} with  $M=O(N^d)$, we have
\begin{equation*}
	E^Q[L^2_{\sigma,\mu,N}(b)] = O(1) N^{2d}\delta_b^{-6}w_{\sigma,\mu,N}^2(b).
\end{equation*}
According to the choice of $N_0$ in Theorem~\ref{bias}, we have
\begin{equation*}
	E^Q[L^2_{\sigma,\mu,N}(b)]=O(1)b^{4d/\beta^*(\frac{3d}{\beta^*}+2+\varepsilon_0)+6}\varepsilon^{-4d/\beta^*-2d\varepsilon_0}
w_{\sigma,\mu,N}^2(b)
\end{equation*}
uniformly for $\mu,\sigma\in \mathcal{C}(\mu_l,\mu_u, \sigma_l,\sigma_u, \beta,\kappa_H)$. This completes our proof.
\end{proof}

\medskip
\begin{proof}[\bf Proof of Corollary~\ref{coro:accu}]
	The mean squared error of $Z_{\sigma,\mu,N}(b)$ is decomposed as the sum of its bias and variance,
\begin{equation*}
	E[Z_{\sigma,\mu,N}(b)-w_{\sigma,\mu}(b)]^2 = [E Z_{\sigma,\mu,N}(b) -w_{\sigma,\mu}(b)]^2 + Var(Z_{\sigma,\mu,N}(b))
	= [ w_{\sigma,\mu,N}(b)- w_{\sigma,\mu}(b)]^2+ Var(L_{\sigma,\mu,N}(b))/n.
\end{equation*}
Setting $\varepsilon:=\varepsilon\delta^{1/2}$ in Theorem~\ref{bias}, we have
$[ w_{\sigma,\mu,N}(b)- w_{\sigma,\mu}(b)]^2<{\varepsilon^2\delta}w_{\sigma,\mu}^2(b)/2$
for $N\geq N(\varepsilon\delta^{1/2},b).
$
Furthermore, according to Theorem~\ref{variance}, we have
$
Var(L_{\sigma,\mu,N}(b))/n\leq {\varepsilon^2\delta}w_{\sigma,\mu}^2(b)/2
$
for $n\geq 2\kappa_c b^{q}\varepsilon^{-q_1-2}\delta^{-\frac{q_1}{2}-1}$.
Consequently, for such $N$ and $n$ we have
$E[Z_{\sigma,\mu,N}(b)-w_{\sigma,\mu}(b)]^2 \leq \varepsilon^2\delta.
$
Thanks to Chebyshev's inequality, we have
\begin{equation*}
	P(|Z_{\sigma,\mu,N}(b)-w_{\sigma,\mu}(b)|>\varepsilon)<\frac{E[Z_{\sigma,\mu,N}(b)-w_{\sigma,\mu}(b)]^2 }{\varepsilon^2}\leq \delta.
\end{equation*}
Therefore, $Z_{\sigma,\mu,N}(b)$ satisfies \eqref{eq:accuracy}.
\end{proof}

\subsection{Proofs of supporting lemmas}\label{sec:proof-lemmas}
\begin{proof}[\bf Proof of Lemma~\ref{lemma:tail-approx}]
	First, according to Lemma~\ref{LemBorell}, we have
	\begin{equation}\label{eq:upper-approx}
		P\Big(\sup_{t\in T}\sigma(t)f(t)+\mu(t)>b\Big )\leq P\left(\sup_{t\in T}\sigma(t)f(t)>b-\max_{t\in T}\mu(t) \right)\leq e^{ -(1+o(1))\frac{b^2}{2\max_{t\in T}\sigma^2(t)}}.
	\end{equation}
	On the other hand, for each $t\in T$ we have
	\begin{equation*}
	P\Big(\sup_{t\in T}\sigma(t)f(t)+\mu(t)>b \Big)\geq P\Big(\sigma(t)f(t)+\mu(t)>b \Big)=P\Big(f(t)>\frac{b-\mu(t)}{\sigma(t)}\Big),
	\end{equation*}
	which is further bounded from below by
	\begin{equation*}
		P\Big(\sup_{t\in T}\sigma(t)f(t)+\mu(t)>b\Big )\geq \frac{1}{\sqrt{2\pi}\sigma(t)}\Big(\frac{\sigma(t)}{b-\mu(t)}-\frac{\sigma^3(t)}{(b-\mu(t))^3}\Big)e^{-\frac{(b-\mu(t))^2}{2\sigma^2(t)}}
		= \tep b^{-1}e^{-\frac{b^2}{2\sigma^2(t)}}.
	\end{equation*}
	To obtain the last equation in the above display, we used the fact that $\mu(t)\in[\mu_l,\mu_u]$ and $\sigma(t)\in[\sigma_l,\sigma_u]$ with $\sigma_l>0$.
	Taking the maximum of the right-hand side of the above display, we have
	\begin{equation}
		P\Big(\sup_{t\in T}\sigma(t)f(t)+\mu(t)>b \Big)\geq \tep b^{-1}\max_{t\in T} e^{-\frac{(b-\mu(t))^2}{2\sigma^2(t)}}.
	\end{equation}
	Combining the above expression with \eqref{eq:upper-approx}, we complete the proof.
\end{proof}
\begin{proof}[\bf Proof of Lemma~\ref{lemma:esup}]
	To prove this lemma, we will need the following entropy bound (\cite{dudley1973}).
	\begin{lemma}\label{lemma:dudley}
Let $f$ be a centered Gaussian field living on a metric space $\mathcal{U}$. Define the pseudo-metric
$$
d_f(s,t) = \sqrt{E (f(s)-f(t))^2}.
$$
Assume that $\mathcal{U}$ is a compact space under the metric $d_f$ and for each $\varepsilon>0$. Denote by  $N(\varepsilon)$ the smallest number of  balls with radius $\varepsilon$ under the metric $d_f$. Then there exists a universal constant $K$ such that
\begin{equation}{}
	E \Big[
	\sup_{t\in\mathcal{U}} f(t)
	\Big]\leq K\int_0^{diam(\mathcal{U})} (\log N(\varepsilon))^{\frac{1}{2}}d\varepsilon.
\end{equation}
	\end{lemma}
Let $\mathcal{U}=\{(s,t):s,t\in T, |s-t|\leq \kappa_m \frac{1}{N}\}$ and   
\begin{equation*}
	d_{\xi}^2(
	(s,t),(s',t')
	)=E[
	\xi(s,t)-\xi(s',t')
	]^2= E[\xi(s,s')-\xi(t,t')]^2.
\end{equation*}
We first investigate the metric $d_{\xi}$. We have
\begin{equation}
	d_{\xi}^2(
	(s,t),(s',t')
	)
	\leq 2Var(\xi(s,s'))+ 2Var(\xi(t,t')).
\end{equation}
Applying \eqref{eq:var-xi} to the above display, we  have that there is a  $\tkappa$ uniformly for all $\sigma,\mu$ satisfying Assumption C1, such that
\begin{equation}\label{eq:d-bound}
	d_{\xi}(
	(s,t),(s',t')
	)
	\leq \tkappa \sqrt{|s-s'|^{\beta^{*}}+|t-t'|^{\beta^{*}}}.
\end{equation}
According to the relationship between the $l_p$ norms, we have
$	(|s-s'|^{\beta^{*}}+|t-t'|^{\beta^{*}})^{\frac{1}{\beta^{*}}}\leq d^{\frac{1}{2}-\frac{1}{\beta^*}}\sqrt{|s-s'|^{2}+|t-t'|^{2}}.
$
The result, together with \eqref{eq:d-bound}, implies that $B((s,t), \tep \varepsilon^{\frac{2}{\beta^{*}}} )\subset B_{d_{\xi}}((s,t),\varepsilon )$ for some constant $\tep$ that only depends on $d, \beta^{*}$ and $\tkappa$, where $B$ and $B_{\xi}$ denote balls under the Euclidean norm and $d_{\xi}$ metrics respectively. Note that the set $T\times T$ can be covered by $\tkappa \varepsilon^{-\frac{4d}{\beta^{*}}}$ many  $B(\tep \varepsilon^{\frac{2}{\beta^{*}}} )$ balls with a possibly different $\tkappa$. Consequently, the set $\mathcal{U}$ can be covered by the same number of $B_{d_{\xi}}(\varepsilon )$ balls. Therefore, we have 
\begin{equation*}
	\log(N(\varepsilon))\leq \log \tkappa + \frac{4d}{\beta^{*}} \log \varepsilon^{-1}
\end{equation*}
On the other hand, we have
$
d_{\xi}(
	(s,t),(s',t')
	)
	\leq 2Var(\xi(s,t))+ 2Var(\xi(s',t')).
$
Also according to \eqref{eq:var-xi}, we have
$	d_{\xi}^2(
	(s,t),(s',t')
	)
	= O(|s-t|^{\beta^*}+|s'-t'|^{\beta^*}).
$
Therefore, for $|s-t|\leq \kappa_m/N$ we have
$	d_{\xi}(
	(s,t),(s',t')
	)=O(N^{-\beta^*/2}).
$ Consequently,
$	diam(\mathcal{U}) \leq \tkappa N^{-\beta^*/2}.
$ According to Lemma~\ref{lemma:dudley}, we have
\begin{equation*}
	E {\sup_{t,s\in T, |t-s|\leq \kappa_m \frac{1}{N}}} \leq  \tkappa (\frac{4d}{\beta^{*}})^{1/2}\int_0^{\tkappa N^{-\beta^*/2}}  (\log \varepsilon^{-1})^{1/2} d\varepsilon = O(N^{-\beta^*/2} \log N).
\end{equation*}
This completes our proof.
\end{proof}

\section*{Acknowledgement}
The authors thank the Editor, an associate Editor and two referees for their constructive comments.
Li's research is partially supported by the National Science Foundation grant DMS-1712657. 
Xu's research is partially supported by the National Science Foundation grants DMS-1712717 and SES-1659328, and National Security Agency grant  
H98230-17-1-0308.

\bibliographystyle{abbrv}
\bibliography{minimalBib,other,bibprob,bibstat}

\end{document}